\documentclass[11pt]{article}
\usepackage[dvips]{graphics}
\usepackage{epsfig,rotating}
\usepackage{amssymb,amsmath}
\usepackage{latexsym}
\usepackage[usenames]{color}

\newcommand{\be}{\begin{equation}}
\newcommand{\ee}{\end{equation}}
\newcommand{\beaa}{\begin{eqnarray*}}
\newcommand{\eeaa}{\end{eqnarray*}}
\newcommand{\bea}{\begin{eqnarray}}
\newcommand{\eea}{\end{eqnarray}}
\newcommand{\lbl}{\label}

\newcommand{\bei}{\begin{itemize}}
\newcommand{\eei}{\end{itemize}}

\newcommand{\bd}{\bold}
\newcommand{\hf}{{1 \over 2}}

\def\bX{\mathbf{X}}
\def\bZ{\mathbf{Z}}

\newtheorem{theorem}{ \noindent T{\footnotesize HEOREM}}
\newtheorem{prop}{ \noindent P{\footnotesize ROPOSITION}}
\newtheorem{lemma}{ \noindent L{\footnotesize EMMA}}[section]
\newtheorem{coro}{ \noindent C{\footnotesize OROLLARY}}[section]
\newtheorem{remark}{ \noindent R{\footnotesize EMARK}}[section]

\newcommand{\RR}{\mathbb{R}}
\newcommand{\goto}{\rightarrow}
\begin{document}

\title{Distributions of Angles in Random Packing on Spheres }
\author{Tony Cai$^{1}$, Jianqing Fan$^2$ and Tiefeng Jiang$^{3}$\\
University of Pennsylvania, Princeton University and University of Minnesota}

\date{}
\maketitle

\footnotetext[1]{Statistics Department, The Wharton School, University of
  Pennsylvania,  Philadelphia, PA 19104, \newline  \indent \ \
tcai@wharton.upenn.edu. The research of Tony Cai was supported in part
by NSF FRG Grant  \newline  \indent \ \
DMS-0854973, NSF Grant DMS-1209166, and NIH Grant R01 CA127334.}
\footnotetext[2]{\noindent Department of Operation Research and Financial Engineering, Princeton University,  Princeton, \newline  \indent \ \ NJ08540,   \ jqfan@princeton.edu. The research of Jianqing Fan was supported in part
by NSF grant \newline  \indent \ \ DMS-1206464  and NIH grants  NIH R01-GM072611 and R01GM100474.}
\footnotetext[3]{\noindent School of Statistics, University of Minnesota, 224 Church
Street, MN55455, tjiang@stat.umn.edu. \newline  \indent \ \
The research of Tiefeng Jiang was
supported in part by NSF FRG Grant DMS-0449365 and  NSF Grant  \newline  \indent \ \
DMS-1209166.}

\begin{abstract}
This paper studies the asymptotic behaviors of the pairwise angles among $n$ randomly and uniformly distributed unit vectors in $\RR^p$ as the number of points $n\goto \infty$, while the dimension $p$ is either fixed or growing with $n$. For both settings, we derive the limiting empirical distribution of  the random angles and the limiting distributions of the extreme angles. The results reveal interesting differences in the two settings and provide a precise characterization of the folklore that ``all high-dimensional random vectors are almost always nearly orthogonal to each other". Applications to statistics and machine learning and connections with some open problems in physics and mathematics are also discussed.

\end{abstract}



\newpage
\section{Introduction}
\lbl{intro}
\setcounter{equation}{0}

The distribution of the Euclidean and geodesic distances between two random points on a unit sphere or other geometric objects has a wide range of applications including transportation networks, pattern recognition, molecular biology, geometric probability, and many branches of physics. The distribution has been well studied in different settings.  For example, Hammersley (1950), Lord (1954), Alagar (1976) and Garc\'{\i}a-Pelayo (2005) studied the distribution of the Euclidean distance between two random points on the unit sphere $\mathbb{S}^{p-1}$.  Williams (2001) showed that,  when the underlying geometric object is a sphere or an ellipsoid, the distribution has a strong connection to the neutron transport theory. Based on applications in neutron star models and tests for random number  generators in $p$-dimensions, Tu and Fischbach (2002) generalized the results from unit spheres to more complex geometric objects including the ellipsoids and discussed many applications.
In general, the angles, areas and volumes associated with random points, random lines and random planes appear in the studies of stochastic geometry, see, e.g., Stoyan and Kendall (2008) and Kendall and Molchanov (2010).

In this paper we consider the empirical law and extreme laws of the pairwise angles among a large number of random unit vectors. More specifically, let $\bd{X}_1, \cdots, \bd{X}_n$ be random points independently chosen with the uniform distribution on $\mathbb{S}^{p-1},$ the unit sphere in $\mathbb{R}^{p}.$ The $n$ points $\bd{X}_1, \cdots, \bd{X}_n$ on the sphere naturally generate $n$ unit vectors $\overset{\longrightarrow}{\bd{O}\bd{X}_i}$ for $i=1,2\cdots,n,$ where $\bd{O}$ is the origin. Let $0\leq \Theta_{ij}\leq \pi$  denote the angle between $\overset{\longrightarrow}{\bd{O}\bd{X}_i}$ and $\overset{\longrightarrow}{\bd{O}\bd{X}_j}$ for all $1\leq i< j\leq n.$  In the case of a fixed dimension, the global behavior of the angles $\Theta_{ij}$ is captured by its empirical distribution
\bea\lbl{empirical}
 \mu_n=\frac{1}{\binom{n}{2}}\sum_{1\leq i < j\leq n}\delta_{\Theta_{ij}},\ n \geq 2.
 \eea
 When both the number of points $n$ and the dimension $p$ grow, it is more appropriate to consider the normalized empirical distribution
 \bea\lbl{empirical1}
 \mu_{n,p}=\frac{1}{\binom{n}{2}}\sum_{1\leq i < j\leq n}\delta_{\sqrt{p-2}({\frac{\pi}{2}-\Theta_{ij}})},\ n \geq 2,\ p\geq 3.
 \eea
In many applications it is of significant interest to consider the extreme angles $\Theta_{\min}$ and $\Theta_{\max}$ defined by
\bea\lbl{geese}
\Theta_{\min}&=&\min\{\Theta_{ij};\ 1\leq i <  j \leq n\};\lbl{ducks}\ \ \ \ \ \ \ \ \ \\
\Theta_{\max}&=&\max\{\Theta_{ij};\ 1\leq i <  j \leq n\}.\lbl{fallway}\ \ \ \ \ \ \ \ \ \
\eea
We will study both the empirical distribution of  the angles $\Theta_{ij}$, $1\le i < j \le n$, and the distributions of the extreme angles $\Theta_{\min}$ and $\Theta_{\max}$ as the number of points $n\goto \infty$, while the dimension $p$ is either fixed or growing with $n$.

The distribution of minimum angle of $n$ points randomly distributed on the $p$-dimensional unit sphere has important implications in statistics and machine learning.  It indicates how strong spurious correlations can be for $p$ observations of $n$-dimensional variables (Fan et al, 2012).  It can be directly used to test isotropic of the distributions (see Section 4).  It is also related to regularity conditions such as the Incoherent Condition (Donoho and Huo, 2001),  the Restricted Eigenvalue Condition (Bickel et al, 2009), the $\ell_q$-Sensitivity (Gautier and Tsybakov, 2011) that are needed for sparse recovery.  See also Section 5.1.

The present paper systematically investigates the asymptotic behaviors of the random angles $\{\Theta_{ij};1\leq i < j \leq n\}$. It is shown that, when the dimension $p$ is fixed, as $n\goto \infty$, the empirical distribution $\mu_n$ converges to a distribution with the density function given by
\[
h(\theta)=\frac{1}{\sqrt{\pi}}\frac{\Gamma(\frac{p}{2})}{\Gamma(\frac{p-1}{2})}\cdot (\sin \theta)^{p-2},\ \theta \in [0, \pi].
\]
On the other hand, when the dimension $p$ grows with $n$,  it is shown that  the limiting normalized empirical distribution $\mu_{n,p}$ of the random angles $\Theta_{ij}$, $1\le i < j \le n$ is Gaussian.  When the dimension is high, most of the  angles are concentrated around $\pi/2$. The results provide a precise description of this concentration and thus give a rigorous theoretical justification to the folklore that ``all high-dimensional random vectors are almost always nearly orthogonal to each other," see, e.g.,   Diaconis and Freedman (1984) and Hall et al (2005). A more precise description is given in Proposition \ref{rabbit} later in terms of the concentration rate.

In addition to the empirical law of the angles $\Theta_{ij}$, we also consider the extreme laws of the random angles in both the fixed  and growing dimension settings. The limiting distributions of  the extremal statistics $\Theta_{\max}$ and $\Theta_{\min}$ are derived. Furthermore, the limiting distribution of the sum of the two extreme angles $\Theta_{\min} + \Theta_{\max}$ is also established. It shows that $\Theta_{\min} + \Theta_{\max}$ is highly concentrated at $\pi$.

The distributions of the minimum and maximum angles as well as the empirical distributions of all pairwise angles have important applications in statistics.  First of all, they can be used to test whether a collection of random data points in the $p$-dimensional Euclidean space follow a spherically symmetric distribution (Fang et al, 1990).  The natural test statistics are either $\mu_n$ or $\Theta_{\min}$ defined respectively in \eqref{empirical} and \eqref{geese}.  The statistic $\Theta_{\min}$ also measures the maximum spurious correlation among $n$ data points in the $p$-dimensional Euclidean space.
The correlations between a response vector with $n$ other variables, based on $n$ observations, are considered as spurious when they are smaller than a certain upper quantile of the distribution of $|\cos(\Theta_{\min})|$ (Fan and Lv, 2008).    The statistic $\Theta_{\min}$ is also related to the bias of estimating the residual variance (Fan et al, 2012).  More detailed discussion of the statistical applications of our studies is given in Section \ref{appl}.

The study of the empirical law and the extreme laws of the random angles $\Theta_{ij}$ is closely connected to several deterministic open problems in physics and mathematics, including the general problem in physics of finding the minimum energy configuration of a system of particles on the surface of a sphere and the mathematical problem of uniformly distributing points on a sphere, which originally arises in complexity theory.  The extreme laws of the random angles   considered in this paper is also related to the study of the coherence of a random matrix, which is defined to be the largest magnitude of the Pearson correlation coefficients between the columns of the random matrix. See Cai and Jiang (2011, 2012) for the recent results and references on the distribution of the coherence. Some of these connections are discussed in more details in Section \ref{discussion.sec}.

This paper is organized as follows. Section \ref{fixedn} studies  the limiting empirical and extreme laws of the angles $\Theta_{ij}$ in the setting of the fixed dimension $p$ as the number of points $n$ going to $\infty.$ The case of growing dimension is considered in Section \ref{largen}. Their applications in statistics are outlined in Section \ref{appl}.   Discussions on the connections to the machine learning and some open problems in physics and mathematics are given in Section \ref{discussion.sec}.
The proofs of the main results  are relegated in Section \ref{proof.sec}.

\section{When The Dimension $p$ Is Fixed}
\lbl{fixedn}


In this section we consider the limiting empirical distribution of the angles $\Theta_{ij}$, $1\le i < j \le n$ when the number of random points $n \goto \infty$ while the dimension $p$ is fixed. The case where both $n$ and $p$ grow will be considered in the next section.
Throughout the paper, we let $\bd{X}_1$, $\bd{X}_2$, $\cdots$, $\bd{X}_n$ be independent random points with the uniform distribution on the unit sphere $\mathbb{S}^{p-1}$ for some fixed $p\ge 2$.

We begin with the limiting empirical distribution of the random angles.
\begin{theorem}[Empirical Law for Fixed $p$]
\lbl{quarter}
Let the empirical distribution $\mu_n$ of the angles $\Theta_{ij}$, $1\le i < j \le n$, be defined as in \eqref{quarter}.
Then, as $n\to\infty$, with probability one, $\mu_n$ converges weakly to the distribution with density
 \bea\lbl{reliable}
 h(\theta)=\frac{1}{\sqrt{\pi}}\frac{\Gamma(\frac{p}{2})}{\Gamma(\frac{p-1}{2})}\cdot (\sin \theta)^{p-2},\ \theta \in [0, \pi].
 \eea
\end{theorem}

In fact, $h(\theta)$ is the probability density function of $\Theta_{ij}$ for any $i \ne j$ ($\Theta_{ij}$'s are identically distributed). Due to the dependency of $\Theta_{ij}$'s, some of them are large and some are small. Theorem \ref{quarter} says that the average of these angles asymptotically has the same density as that of  $\Theta_{12}$.

Notice that when $p = 2$, $h(\theta)$ is the uniform density on $[0, \pi]$, and when $p > 2$, $h(\theta)$ is unimodal with mode $\theta=\pi/2$.   Theorem \ref{quarter} implies that most of the angles in the total of $\binom{n}{2}$ angles are concentrated around $\pi/2$.  This concentration becomes stronger as the dimension $p$ grows since $(\sin \theta)^{p-2}$ converges to zero more quickly for $\theta \ne \pi/2$.
In fact, in the extreme case when $p\to \infty$, almost all of  $\binom{n}{2}$  angles go to $\pi/2$ at the rate $\sqrt{p}$.   This can  be seen from Theorem \ref{nickle} later.

\begin{figure}[htp]
\centerline{\includegraphics[scale=1]{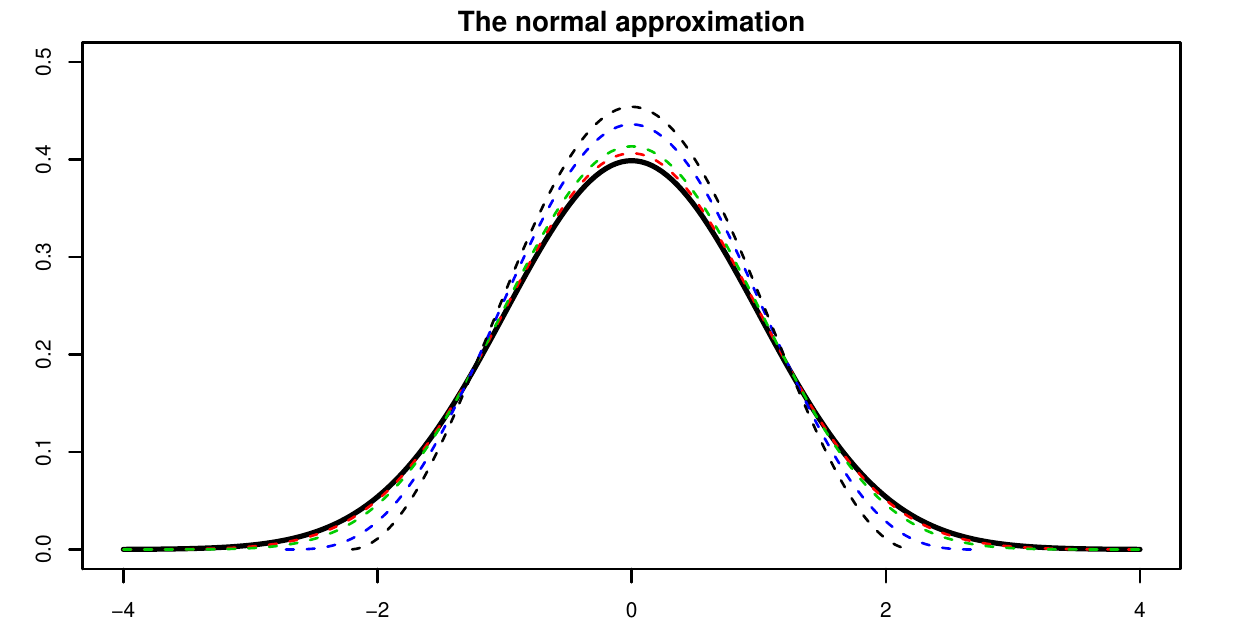}}
\caption{Functions $h_p(\theta)$ given by (\ref{eq1}) for $p=4, 5, 10$ and $20$.  They are getting closer to the normal density (thick black) as $p$ increases.}
\label{fig1}
\end{figure}

It is helpful to see how the density changes with the dimension $p$. Figure~\ref{fig1} plots the function
\begin{eqnarray} \label{eq1}
 h_p(\theta) & = & \frac{1}{\sqrt{p-2}} h\Bigl (\frac{\pi}{2}- \frac{\theta}{\sqrt{p-2}} \Bigr ) \nonumber \\
  & = &
 \frac{1}{\sqrt{\pi}}\frac{\Gamma(\frac{p}{2}) }{\Gamma(\frac{p-1}{2})  \sqrt{p-2}}\cdot \Bigl (\cos \frac{\theta}{\sqrt{p-2}} \Bigr )^{p-2},\quad \theta \in [0, \pi]
\end{eqnarray}
which is the asymptotic density of the normalized empirical distribution $\mu_{n,p}$ defined in \eqref{empirical1} when the dimension $p$ is fixed. Note that in the definition of  $\mu_{n,p}$ in \eqref{empirical1}, if  ``$\sqrt{p-2}$" is replaced by ``$\sqrt{p}$", the limiting  behavior of $\mu_{n,p}$ does not change when both $n$ and $p$ go to inifnity. However,  it shows in our simulations and the approximation (\ref{eq2}) that the fitting is better for  relatively small  $p$ when  ``$\sqrt{p-2}$" is used.

Figure~\ref{fig1} shows that the distributions $h_p(\theta) $ are very close to normal when $p\geq 5$.  This can also be seen from the asymptotic approximation
\begin{equation} \label{eq2}
h_p(\theta) \propto  \exp \Bigl ((p-2) \log\bigl \{ \cos\bigl (\frac{\theta}{\sqrt{p-2}} \bigr) \bigr \} \Bigr) \approx  e^{-\theta^2/2}.
\end{equation}


We now consider the limiting distribution of the extreme angles $\Theta_{\min}$ and $\Theta_{\max}$.

\begin{theorem}[Extreme Law for Fixed $p$]
\lbl{tube}
Let $\Theta_{\min}$ and  $\Theta_{\max}$ be defined as in \eqref{ducks} and \eqref{fallway} respectively.
Then,  both $n^{2/(p-1)}\Theta_{\min}$ and $n^{2/(p-1)}(\pi-\Theta_{\max})$ converge weakly to a distribution given by
\bea\lbl{chick}
F(x)=\begin{cases}
1-e^{-Kx^{p-1}},& \text{ if $x\geq 0$;}\\
0, & \text{ if $x< 0$,}
\end{cases}
\eea
 as $n\to\infty$, where
\bea\lbl{luck}
K=\frac{1}{4\sqrt{\pi}}\frac{\Gamma(\frac{p}{2})}{\Gamma(\frac{p+1}{2})}.
\eea
\end{theorem}
The above theorem says that the smallest angle $\Theta_{\min}$ is close to zero, and the largest angle $\Theta_{\max}$ is close to $\pi$ as $n$ grows. This makes sense from Theorem \ref{quarter} since the support of the density function $h(\theta)$ is $[0, \pi].$

In the special case of  $p=2,$ the scaling of $\Theta_{\min}$ and $\pi-\Theta_{\max}$ in Theorem \ref{tube} is $n^2.$ This is in fact can also be seen in a similar problem. Let $\xi_1, \cdots, \xi_n$ be i.i.d.  $U[0, 1]$-distributed random variables with the order statistics $\xi_{(1)}\le \cdots \le \xi_{(n)}.$ Set $W_n:=\min_{1\leq i \leq n-1}(\xi_{(i+1)}-\xi_{(i)}),$ which is the smallest spacing among the observations of $\xi_i$'s. Then, by using the representation theorem of $\xi_{(i)}$'s through i.i.d. random variables with exponential distribution ${\rm Exp}(1)$ (see, e.g., Proposition 4.1 from Resnick (1987)), it is easy to check that $n^2W_n$ converges weakly to ${\rm Exp}(1)$ with the probability density function $e^{-x}I(x\geq 0).$

\begin{figure}[htbp]
\centerline{\includegraphics[scale=1]{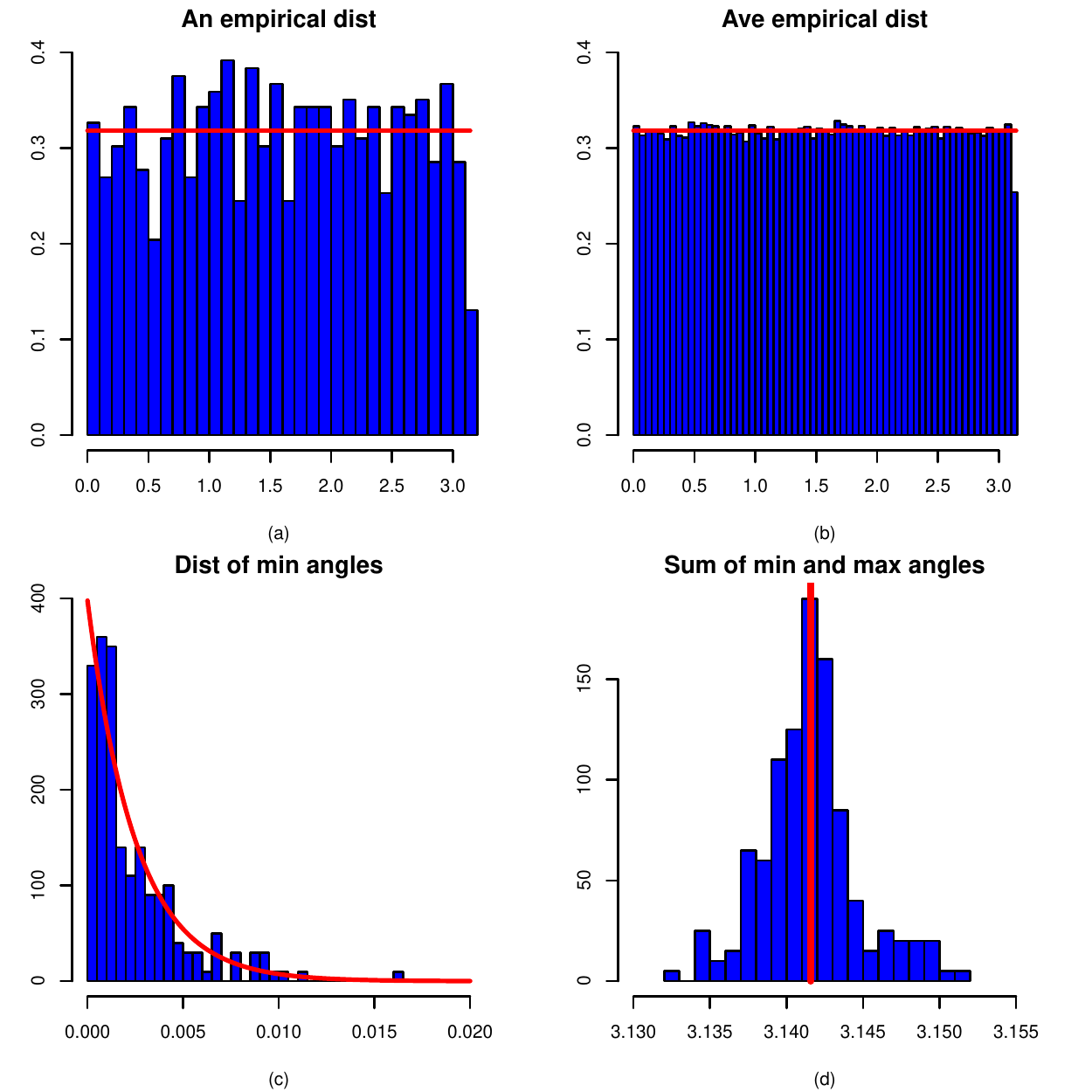}}
\caption{Various distributions for $p=2$ and $n=50$ based on 200 simulations. (a) A realization of the empirical distribution $\mu_n$;  (b) The average distribution of 200 realizations of $\mu_n$; (c) the distribution of $\Theta_{\min}$ and its asymptotic distribution $\exp(-x/(2 \pi))/(2 \pi)$; (d) the distribution of $\Theta_{\min} + \Theta_{\max}$; the vertical line indicating the location $\pi$.}
\label{fig2}
\end{figure}

\begin{figure}[htbp]
\centerline{\includegraphics[scale=1]{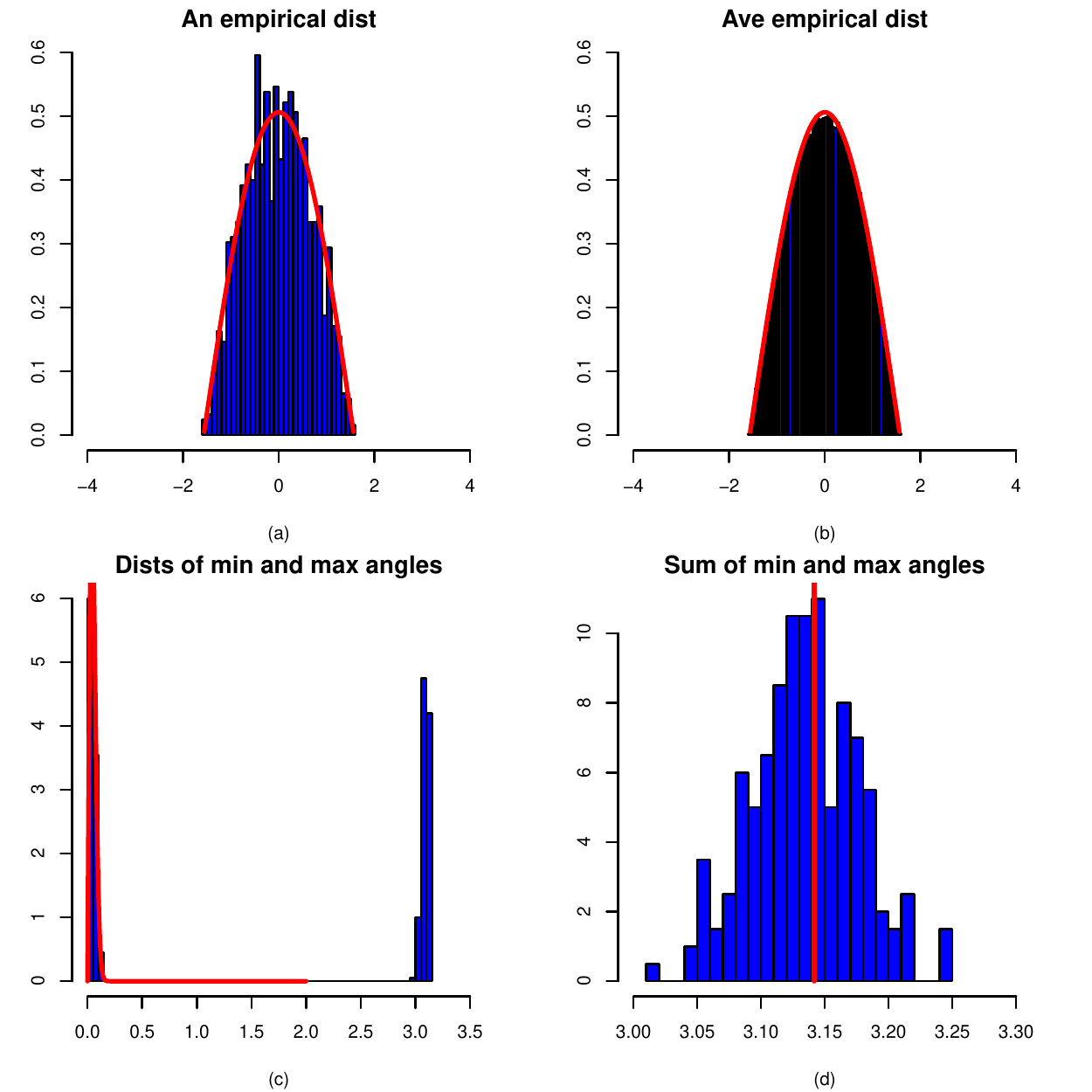}}
\caption{Various distributions for $p=3$ and $n=50$ based on 200 simulations. (a) A realization of the normalized empirical distribution $\mu_{n, p}$ given by (\ref{empirical1});  (b) The average distribution of 200 realizations of $\mu_{n,p}$; (c) the distribution of $\Theta_{\min}$ and its asymptotic distribution; (d) the distribution of $\Theta_{\min} + \Theta_{\max}$; the vertical line indicating the location $\pi$.}
\label{fig3}
\end{figure}

\begin{figure}[htbp]
\centerline{\includegraphics[scale=1]{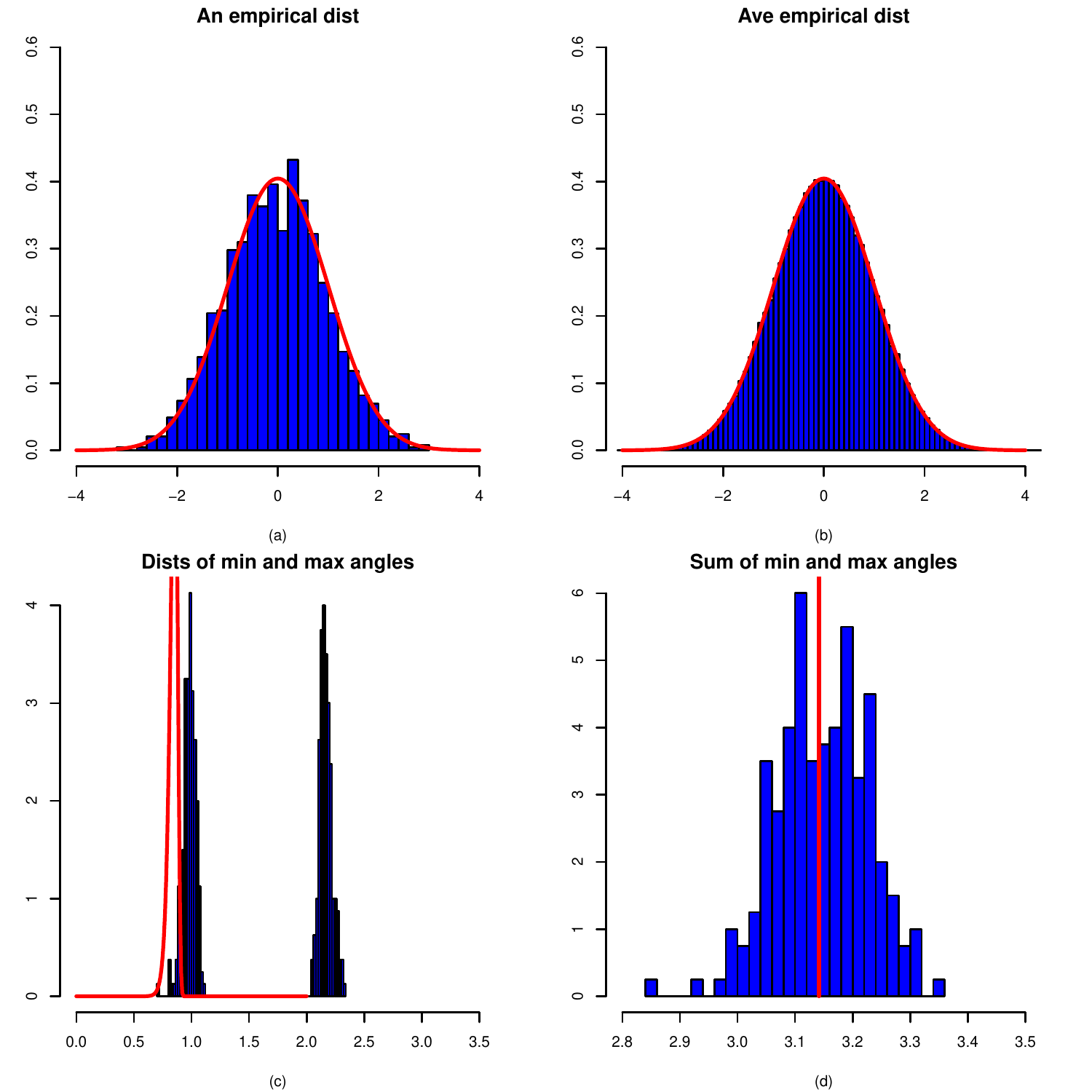}}
\caption{Various distributions for $p=30$ and $n=50$ based on 200 simulations. (a) A realization of the normalized empirical distribution $\mu_{n, p}$ given by (\ref{empirical1});  (b) The average distribution of 200 realizations of $\mu_{n,p}$; (c) the distribution of $\Theta_{\min}$ and its asymptotic distribution; (d) the distribution of $\Theta_{\min} + \Theta_{\max}$; the vertical line indicating the location $\pi$.}
\label{fig4}
\end{figure}

To see the goodness of the finite sample approximations, we simulate 200 times from the distributions with $n=50$ for $p=2, 3$ and 30.  The results are shown respectively in Figures~\ref{fig2}--\ref{fig4}.  Figure 2 depicts the results when $p=2$.  In this case, the empirical distribution $\mu_n$ should approximately be uniformly distributed on $[0, \pi]$ for most of realizations.  Figure~\ref{fig2} (a) shows that it holds approximately truly for $n$ as small as 50 for a particular realization (It indeed holds approximately for almost all realizations).  Figure~\ref{fig2}(b) plots the average of these 200 distributions, which is in fact extremely close to the uniform distribution on $[0, \pi]$. Namely, the bias is negligible.  For $\Theta_{\min}$, according to Theorem~\ref{quarter}, it should be well approximated by an exponential distribution with $K = 1/(2\pi)$.  This is verified by Figure~\ref{fig2}(c), even when sample size is as small as 50.  Figure~\ref{fig2}(d) shows the distribution of $\Theta_{\min} + \Theta_{\max}$ based on the 200 simulations.  The sum is distributed tightly around $\pi$, which is indicated by the red line there.

The results for $p=3$ and $p=30$ are demonstrated in Figures~\ref{fig3} and \ref{fig4}.  In this case, we show the empirical distributions of $\sqrt{p-2} (\pi/2-\Theta_{ij} )$  and their asymptotic distributions.  As  in Figure~\ref{fig1}, they are normalized. Figure~\ref{fig3}(a) shows a realization of the distribution and Figure~\ref{fig3}(b) depicts the average of 200 realizations of these distributions for $p = 3$.  They are very close to the asymptotic distribution, shown in the curve therein.  The distributions of $\Theta_{\min}$ and $\Theta_{\max}$ are plotted in Figure~\ref{fig3}(c).  They concentrate respectively around $0$ and $\pi$.  Figure~\ref{fig3}(d) shows that the sum is concentrated symmetrically around $\pi$.

When $p = 30$, the approximations are still very good for the normalized empirical distributions.  In this case, the limiting distribution is indistinguishable from the normal density, as shown in Figure~\ref{fig1}.  However, the distribution of $\Theta_{\min}$ is not approximated well by its asymptotic counterpart, as shown in Figure~\ref{fig4}(c).  In fact,
$\Theta_{\min}$ does not even tends to zero.  This is not entirely surprising since $p$ is comparable with $n$.  The asymptotic framework in Section~\ref{largen} is more suitable.  Nevertheless, $\Theta_{\min} + \Theta_{\max}$ is still symmetrically distributed around $\pi$.

The simulation results show  that $\Theta_{\max}+\Theta_{\min}$ is very close to $\pi.$ This  actually can be seen trivially from Theorem \ref{tube}: $\Theta_{\min}\to 0$ and $\Theta_{\max}\to \pi$ in probability as $p\to \infty.$ Hence, the sum goes to $\pi$ in probability.  An interesting question is: how fast is this convergence?  The following result  answers this question.

\begin{theorem}[Limit Law for Sum of Largest and Smallest Angles]
\lbl{goose}
Let $\bd{X}_1$, $\bd{X}_2$, $\cdots$, $\bd{X}_n$ be independent random points with the uniform distribution on  $\mathbb{S}^{p-1}$ for some fixed $p\ge 2$. Let $\Theta_{\min}$ and  $\Theta_{\max}$ be defined as in \eqref{ducks} and \eqref{fallway} respectively.
Then, $n^{2/(p-1)}\big(\Theta_{\max}+\Theta_{\min}-\pi\big)$  converges weakly to the distribution of $X-Y$, where $X$ and $Y$ are i.i.d. random variables with distribution function $F(x)$ given in (\ref{chick}).
\end{theorem}

It is interesting to note that the marginal distribution of $\Theta_{\min}$ and
$\pi-\Theta_{\max}$ are identical.  However, $n^{2/(p-1)} \Theta_{\min}$ and  $n^{2/(p-1)}(\pi-\Theta_{\max})$ are asymptotically independent with non-vanishing limits and hence their difference is non-degenerate. Furthermore, since $X$ are $Y$ are i.i.d., $X-Y$ is a symmetric random variable. Theorem \ref{goose} suggests that  $\Theta_{\max}+\Theta_{\min}$ is larger or smaller than $\pi$ ``equally likely".
The symmetry of the distribution of $\Theta_{\max}+\Theta_{\min}$ has already been demonstrated in Figures~\ref{fig2} -- \ref{fig4}.

\section{When Both $n$ and $p$ Grow}
\lbl{largen}

We now turn to the case where  both $n$ and $p$ grow.  The following result shows that the empirical distribution of the random angles, after  suitable normalization, converges to a standard normal distribution. This is clearly different from the limiting distribution given in Theorem \ref{quarter} when the dimension $p$ is fixed.

\begin{theorem}[Empirical Law for Growing $p$]
\lbl{nickle}
Let $\mu_{n,p}$ be defined as in \eqref{empirical1}. Assume $\lim_{n\to\infty}p_n=\infty.$ Then, with probability one, $\mu_{n,p}$ converges weakly to $N(0,1)$ as $n\to\infty.$
\end{theorem}


Theorem \ref{nickle} holds regardless of the speed of $p$ relative to $n$ when both go to infinity.  This has also been empirically demonstrated in Figures~\ref{fig2}--\ref{fig4} (see plots (a) and (b) therein).  The theorem  implies that  most  of  the $\binom{n}{2}$  random angles go to $\pi/2$ very quickly.  Take any $\gamma_p \goto 0$ such that $\sqrt{p} \gamma_p \goto \infty$ and denote by $N_{n,p}$ the number of the angles  $\Theta_{ij}$ that are within $\gamma_p$ of $\pi/2$, i.e., $|\frac{\pi}{2}-\Theta_{ij}|\le \gamma_p$. Then $N_{n,p}/{n\choose 2} \goto 1$. Hence, most of the random vectors in the high-dimensional Euclidean spaces  are nearly orthogonal.  An interesting question is: Given two such random vectors, how fast is their angle close to $\pi/2$ as the dimension increases? The following result answers this question.

\begin{prop}\lbl{rabbit} Let $\bd{U}$ and $\bd{V}$ be two random points on the unit sphere in $\mathbb{R}^p.$ Let $\Theta$ be the angle between $\overset{\longrightarrow}{\bd{O}\bd{U}}$ and $\overset{\longrightarrow}{\bd{O}\bd{V}}.$ Then
\beaa
P(|\Theta -\frac{\pi}{2}|\geq \epsilon) \leq K\sqrt{p}(\cos \epsilon)^{p-2}
\eeaa
for all $p\geq 2$ and $\epsilon \in (0, \pi/2),$  where $K$ is a universal constant.
\end{prop}
Under the  spherical invariance one can think of $\Theta$  as a function of the random point $\bd{U}$ only. There are general concentration inequalities on such functions, see, e.g., Ledoux (2005). Proposition \ref{rabbit} provides a more precise inequality.

 One can see that, as the  dimension $p$ grows, the probability decays exponentially. In particular, take $\epsilon = \sqrt{(c\log p)/p}$ for some constant $c>1$. Note that $\cos \epsilon \le 1- \epsilon^2/2 +\epsilon^4/24$, so
\be\lbl{rainbow}
P\left(|\Theta -\frac{\pi}{2}|\geq \sqrt{\frac{c\log p}{p}}\right)
 \leq K\sqrt{p}\left(1-\frac{c\log p}{2p} + \frac{c^2\log^2 p}{24 p^2} \right)^{p-2}\le K' p^{-\hf(c-1)}
\ee
for all sufficiently large $p$, where $K'$ is a constant depending only on $c$. Hence, in the high dimensional space, the angle between two random vectors is within
$\sqrt{(c\log p)/p}$ of $\pi/2$ with high probability. This provides a precise characterization of the folklore mentioned earlier that ``all high-dimensional random vectors are almost always nearly orthogonal to each other".

We now turn to the limiting extreme laws of the angles when both $n$ and $p \goto \infty$.  For the extreme laws, it is necessary to  divide into three asymptotic regimes: sub-exponential case $\frac{1}{p}\log n \to 0$, exponential case $\frac{1}{p}\log n \to \beta\in (0, \infty)$, and super-exponential case $\frac{1}{p}\log n  \to\infty$. The limiting extreme laws are different in these three regimes.

\begin{theorem}[Extreme Law: Sub-Exponential Case]\lbl{newkind}
Let $p=p_n \to \infty$ satisfy $\frac{\log n}{p} \to 0$ as $n\to\infty$.  Then
\bei
\item[{\rm (i).}] $\max_{1\leq i < j\leq n}|\Theta_{ij} - \frac{\pi}{2}|\to 0$ in probability as $n\to\infty;$

\item[{\rm (ii).}] As $n\to\infty$, $2p\log \sin \Theta_{\min} + 4\log n -\log \log n$
converges weakly to the extreme value distribution  with the distribution function $F(y)=1- e^{-Ke^{y/2}},\  y\in\mathbb{R}$ and $K=1/(4\sqrt{2\pi}\,).$ The conclusion still holds if $\Theta_{\min}$ is replaced by $\Theta_{\max}$.
\eei
\end{theorem}
In this case, both $\Theta_{\min}$ and $\Theta_{\max}$ converge to $\pi/2$ in probability.
The above extreme value distribution differs from that in (\ref{chick}) where the dimension $p$ is fixed. This is obviously caused by the fact that $p$ is finite in Theorem \ref{tube} and goes to infinity in Theorem \ref{newkind}.
\begin{coro}
\lbl{newbrother}
Let $p=p_n$ satisfy $\lim_{n\to\infty}\frac{\log n}{\sqrt{p}}=\alpha\in [0, \infty)$. Then $p\cos^2\Theta_{\min} - 4\log n +\log \log n$ converges weakly to a distribution with the cumulative distribution function $\exp\{-{1\over 4\sqrt{2\pi}} e^{-(y+8\alpha^2)/2}\}$,  $y\in \mathbb{R}$. The conclusion still holds if $\Theta_{\min}$ is replaced by $\Theta_{\max}$.
\end{coro}

\begin{theorem}[Extreme Law: Exponential Case]\lbl{newjean}
Let $p=p_n$ satisfy $\frac{\log n}{p} \to \beta \in (0, \infty)$ as $n\to\infty$, then
\bei
\item[{\rm (i).}] $\Theta_{\min}\to \cos^{-1}\sqrt{1-e^{-4\beta}}$ and $\Theta_{\max}\to \pi-\cos^{-1}\sqrt{1-e^{-4\beta}}$ in probability as $n\to\infty;$

\item[{\rm (ii).}] As $n\to\infty$, $2p\log \sin \Theta_{\min} + 4\log n -\log \log n$ converges weakly to a distribution  with the distribution function
\bea\lbl{neweast}
F(y)=1- \exp\left\{-K(\beta) e^{(y+8\beta)/2}\right\},\  y\in\mathbb{R},\ \mbox{where}\    K(\beta)=\Big(\frac{\beta}{8\pi(1- e^{-4\beta})}\Big)^{1/2},
\eea
 and the conclusion still holds if $\Theta_{\min}$ is replaced by $\Theta_{\max}$.
\eei
\end{theorem}

In contrast to Theorem \ref{newkind}, neither $\Theta_{\max} $ nor $\Theta_{\min}$ converges to $\pi/2$ under the case that $(\log n)/p \to \beta\in (0, \infty).$ Instead, they converge to different constants depending on $\beta.$
\begin{theorem}[Extreme Law: Super-Exponential Case]
\lbl{newlearn}
Let $p=p_n$ satisfy $\frac{\log n}{p} \to \infty$ as $n\to\infty$. Then,
\bei
\item[{\rm (i).}] $\Theta_{\min}\to 0$ and  $\Theta_{\max}\to \pi$ in probability as $n\to\infty;$

\item[{\rm (ii).}] As $n\to\infty$, $2p\log \sin \Theta_{\min} + \frac{4p}{p-1}\log n-\log p$
 converges weakly to the extreme value distribution with the distribution function $F(y)=1- e^{-Ke^{y/2}},\  y\in\mathbb{R}$ with $K=1/(2\sqrt{2\pi}).$ The conclusion still holds if $\Theta_{\min}$ is replaced by $\Theta_{\max}$.
\eei
\end{theorem}

It can be seen from Theorems \ref{newkind}, \ref{newjean} and \ref{newlearn} that $\Theta_{\max}$ becomes larger when the rate $\beta=\lim(\log n)/p$ increases. They are $\pi/2$, $\pi-\cos^{-1}\sqrt{1-e^{-4\beta}}\in (\pi/2, \pi)$ and $\pi$ when $\beta=0$, $\beta\in (0, \infty)$ and $\beta=\infty,$ respectively.

Set $f(\beta)=\pi-\cos^{-1}\sqrt{1-e^{-4\beta}}.$ Then $f(0)=\pi/2$ and $f(+\infty)=\pi$, which corresponds to $\Theta_{\max}$ in (i) of Theorem \ref{newkind} and  (i) of Theorem \ref{newlearn}, respectively. So the conclusions in Theorems \ref{newkind}, \ref{newjean} and \ref{newlearn} are consistent.

Theorem \ref{goose} provides the limiting distribution of $\Theta_{\max}+\Theta_{\min}-\pi$ when the dimension $p$ is fixed. It is easy to see from the above theorems that  $\Theta_{\max}+\Theta_{\min}-\pi\to 0$ in probability as both $n$ and $p$ go to infinity. Its asymptotic distribution is much more involved and we leave it as future work.

\begin{remark}\lbl{fat_man}{\rm
As mentioned in the introduction, Cai and Jiang (2011, 2012) considered the limiting distribution of the coherence of a random matrix and the coherence is closely related to the minimum angle $\Theta_{\min}$. In the current setting, the coherence $L_{n,p}$ is defined by
\be\lbl{mount3}
L_{n,p}=\max_{1\leq i <  j \leq n}|{\rho}_{ij}|
\ee
where ${\rho}_{ij}=\bd{X}_i^T\bd{X}_j$. The results in Theorems \ref{newkind}, \ref{newjean} and \ref{newlearn} are new. Their proofs can be essentially reduced to the analysis  of $\max_{1\leq i < j\leq n}\rho_{ij}$. This maximum is analyzed through modifying the proofs  of the results for the limiting distribution of the coherence $L_{n,p}$ in  Cai and Jiang (2012). The key step in the proofs is the study of the maximum and minimum of {\it pairwise} i.i.d. random variables $\{\rho_{ij};\ 1\leq i < j\leq n\}$ by using the Chen-Stein method.  It is noted that $\{\rho_{ij};\ 1\leq i < j\leq n\}$ are not i.i.d. random variables (see, e.g.,  p.148 from Muirhead (1982)), the standard techniques to analyze the extreme values of  $\{\rho_{ij};\ 1\leq i < j\leq n\}$  do not apply.

}
\end{remark}

\section{Applications to Statistics}
\label{appl}

The results developed in the last two sections can be applied  to test the spherical symmetry (Fang et al, 1990):
\begin{equation} \label{d1}
  H_0:  \mbox{$\bZ$ is spherically symmetric in $\RR^p$}
\end{equation}
based on an i.i.d. sample $\{\bZ_i\}_{i=1}^n$.
Under the null hypothesis  $H_0$,  $\bZ/\|\bZ\|$ is uniformly distributed on $\mathbb{S}^{p-1}$.  It is expected that the minimum angle $\Theta_{\min}$ is stochastically larger under the null hypothesis than that under the alternative hypothesis. Therefore, one should reject the null hypothesis when $\Theta_{\min}$ is too small or formally, reject $H_0$ when
\begin{equation} \label{d2}
  n^{2/(p-1)} \Theta_{\min} \leq c_\alpha,
\end{equation}
where the critical value $c_\alpha$, according to Theorem 2, is given by
$$
   c_\alpha = \left (- K^{-1} \log(1-\alpha) \right)^{1/(p-1)}
$$
for the given significance level $\alpha$.  This provides the minimum angle test for sphericity or the packing test on sphericity.

We run a simulation study to examine the power of the packing test.  The following 6 data generating processes are used:
\begin{enumerate} \itemindent -0.2 in
\item []{\bf Distribution 0}:  the components of $\bX$ follow independently the standard normal distribution;
\item [] {\bf Distribution 1}:  the components of $\bX$ follow independently the uniform distribution on $[-1, 1]$;
\item [] {\bf Distribution 2}:  the components of $\bX$ follow independently the uniform distribution on $[0, 1]$;
\item [] {\bf Distribution 3}:  the components of $\bX$ follow the standard normal distribution with correlation 0.5;
\item [] {\bf Distribution 4}:  the components of $\bX$ follow the standard normal distribution with correlation 0.9;
\item [] {\bf Distribution 5}:  the components of $\bX$ follow independently the mixture distribution $2/3\exp(-x) I(x \geq 0) + 1/3 \exp(x) I(x \leq 0)$.
\end{enumerate}
The results are summarized in Table~\ref{tab1} below.  Note that for Distribution 0, the power corresponds to the size of the test, which is slightly below $\alpha = 5\%$.

\begin{table}[ht]
\begin{center}
\caption{The power (percent of rejections) of the packing test based on 2000 simulations}
\begin{tabular}{crrrrrr}
\hline
Distribution &  0    &   1  &   2    &   3   &   4   &   5   \\ \hline
$p = 2$      &  4.20 & 5.20 &  20.30 & 5.55  & 10.75 &  5.95  \\
$p = 3$      &  4.20 & 6.80 &  37.20 & 8.00  & 30.70 &  8.05  \\
$p = 4$      &  4.80 & 7.05 &  64.90 &11.05  & 76.25 & 11.20 \\
$p = 5$      &  4.30 & 7.45 &  90.50 &18.25  & 99.45 & 11.65 \\
\hline
\end{tabular}
\label{tab1}
\end{center}
\end{table}

The packing test does not examine whether there is a gap in the data on the sphere.
An alternative test statistic is $\mu_n$ or its normalized version $\mu_{n, p}$ when $p$ is large, defined respectively by (\ref{empirical}) and (\ref{empirical1}). A natural test statistic is then to use a distance such as the Kolmogrov-Smirnov distance between $\mu_n$ and $h(\theta)$.  In this case, one needs to derive further the null distribution of such a test statistic.  This is beyond the scope of this paper and we leave it for future work.

Our study also shed lights on the magnitude of spurious correlation.  Suppose that we have a response variable $Y$ and its associate covariates $\{X_j\}_{j=1}^p$ (e.g., gene expressions).  Even when there is no association between the response and the covariate, the maximum sample correlation between $X_j$ and $Y$ based on a random sample of size $n$ will not be zero.  It is closely related to the minimum angle $\Theta_{\min}$ (Fan and Lv, 2008).  Any correlation below a certain thresholding level can be spurious -- the correlation of such a level can occur purely by chance.  For example, by Theorem 6(ii), any correlation (in absolute value) below
\begin{equation}   \label{d4}
   \sqrt{1 - n^{-4/p} (\log (n))^{1/p}}
\end{equation}
can be regarded as the spurious one.  Take, for example,  $p = 30$ and $n = 50$ as in Figure~\ref{fig4}, the spurious correlation can be as large 0.615 in this case.

The spurious correlation also helps understand the bias in calculating the residual $\sigma^2 = \mbox{var}(\varepsilon)$ in the sparse linear model
\begin{equation} \label{d5}
     Y = \mathbf{X}_S^T
     \mathbf{\beta}_S + \varepsilon
\end{equation}
where $S$ is a subset of variables $\{1, \cdots p\}$.  When an extra variable besides $X_S$ is recruited by a variable selection algorithm, that extra variable is recruited to best predict $\varepsilon$ (Fan et al, 2012).   Therefore, by the classical  formula for the residual variance, $\sigma^2$ is underestimated by a factor of $1-\cos^2(\Theta_{\min})$.  Our asymptotic result gives the order of magnitude of such a bias.

\section{Discussions}
\label{discussion.sec}

We have established the limiting empirical and extreme laws of the angles between random unit vectors, both for the fixed dimension and growing dimension cases.  For fixed $p$, we study the empirical law of angles, the extreme law of angles and the law of the sum of the largest and smallest angles in Theorems \ref{quarter}, \ref{tube} and \ref{goose}. Assuming $p$ is large, we establish the empirical law of random angles in Theorem \ref{nickle}. Given two vectors $\bd{u}$ and $\bd{v}$, the cosine of their angle is equal to the Pearson correlation coefficient between them. Based on this observation, among the results developed in this paper, the limiting distribution of the minimum angle $\Theta_{\min}$ given in Theorems \ref{newkind}-\ref{newlearn} for the setting where both $n$ and $p \goto \infty$ is obtained by similar arguments to those in Cai and Jiang (2012) on the coherence of an $n\times p$ random matrix (a detailed discussion is given in Remark \ref{fat_man}). See also Jiang (2004), Li and Rosalsky (2006), Zhou (2007), Liu et al (2008), Li et al (2009) and Li et al (2010) for earlier results on the distribution of the coherence which were all established under the assumption that both $n$ and $p \goto \infty$.

The study of the random angles $\Theta_{ij}$'s,  $\Theta_{\min}$ and $\Theta_{\max}$ is also related to several problems in machine learning as well as some deterministic open problems in physics and mathematics. We briefly discuss some of these connections below.

\subsection{Connections to Machine Learning}
Our studies shed lights on random geometric graphs, which are formed by $n$ random points on the $p$-dimensional unit sphere as vertices with edge connecting between points $\bX_i$ and $\bX_j$ if $\Theta_{ij} > \delta$ for certain $\delta$ (Penrose, 2003; Devroy et al, 2011).  Like testing isotropicity in Section 4, a generalization of our results can be used to detect if there are any implanted cliques in a random graph, which is a challenging problem in machine learning.  It can also be used to describe the distributions of the number of edges and degree of such a random geometric graph.  Problems of hypothesis testing on isotropicity of covariance matrices have strong connections with clique numbers of geometric random graphs as demonstrated in the recent manuscript by 
Castro et al (2012). This furthers connections of our studies in Section 4 to this machine learning problem.

Principal component analysis  (PCA) is one of the most important techniques in high-dimensional data analysis for visualization, feature extraction, and dimension reduction. It has a wide range of applications in statistics and machine learning. A key aspect of the study of PCA in the high-dimensional setting is the understanding of the properties of the principal  eigenvectors of the sample covariance matrix. In a recent paper, Shen et al (2013) showed an interesting asymptotic conical structure in the critical sample eigenvectors under a spike covariance models when the ratio between the dimension and the product of the sample size with the spike size converges to a nonzero constant. They showed that in such a setting the critical sample eigenvectors lie in a right circular cone around the corresponding population eigenvectors. Although these sample eigenvectors converge to the cone, their locations within the cone are random. The behavior of the randomness of the eigenvectors within the cones is related to the behavior of the random angles studied in the present paper.  It is of significant interest to rigorously explore these connections. See Shen et al (2013)  for further discussions.

\subsection{Connections to Some Open Problems in Mathematics and Physics }

The results on random angles established in this paper can be potentially used to study a number of open deterministic problems in mathematics and physics. 

Let $\bd{x}_1, \cdots,\bd{x}_n$ be $n$ points on $\mathbb{S}^{p-1}$ and $R=\{\bd{x}_1, \cdots,\bd{x}_n\}.$ The $\alpha$-energy function is defined by
\beaa\lbl{energy}
E(R, \alpha)
=\begin{cases}
 \sum_{1\leq i < j \leq n}\|\bd{x}_i-\bd{x}_j\|^{\alpha},& \ \text{if $\alpha \ne 0;$}\\
 \sum_{1\leq i < j \leq n}\log \frac{1}{\|\bd{x}_i-\bd{x}_j\|},& \  \text{if $\alpha = 0$,}
\end{cases}
\eeaa
and
$
E(R, -\infty)=\min_{1\leq i < j \leq n}\frac{1}{\|\bd{x}_i-\bd{x}_j\|}
$
where $\|\cdot\|$ is the Euclidean norm in $\mathbb{R}^p.$
 These are known as the electron problem ($\alpha = 0$) and the Coulomb potential problem ($\alpha = -1$). See, e.g., Kuijlaars and Saff (1998) and Katanforoush and
Shahshahani (2003).
The goal is to find the extremal $\alpha$-energy
\beaa\lbl{bow}
\epsilon(R, \alpha):=
\begin{cases}
\inf_{R} E(R, \alpha), & \ \text{if $\alpha \leq 0,$}\\
\sup_{R} E(R, \alpha),  & \ \text{if $\alpha > 0,$}
\end{cases}
\eeaa
and the extremal configuration $R$ that attains $\epsilon(R, \alpha)$.
In particular, when $\alpha=-1,$ the quantity $\epsilon(R, -1)$ is the minimum of the Coulomb potential
\beaa\lbl{sign}
\sum_{1\leq i < j \leq n} \frac{1}{\|\bd{x}_i-\bd{x}_j\|}.
\eeaa
These open problems, as a function of $\alpha,$ are: (i) $\alpha= -\infty$: Tammes problem;
(ii) $\alpha= -1$: Thomson problem;
(iii) $\alpha= 1$: maximum average distance problem;
and (iv) $\alpha= 0$: maximal product of distances between all pairs.
Problem (iv) is the 7th of the 17 most challenging mathematics problems in the 21st century according to Smale (2000). See, e.g., Kuijlaars and Saff (1998) and  Katanforoush  and Shahshahani (2003), for further details.

The above problems can also be formulated through randomization.  Suppose that $\bd{X}_1, \cdots, \bd{X}_n$ are i.i.d. uniform random vectors on $\mathbb{S}^{p-1}$. Suppose $R=\{\bd{x}_1, \cdots,\bd{x}_n\}$ achieves the infinimum or supremum in  the definition of $\epsilon(R, \alpha)$. Since $P(\max_{1\leq i \leq n}\|\bd{X}_i-\bd{x}_i\|<\epsilon)>0$ for any $\epsilon>0,$ it is easy to see that $\epsilon(R, \alpha)=\mbox{ess}\cdot \inf(E(R, \alpha))$ for $\alpha \leq 0$ and  $\epsilon(R, \alpha)=\mbox{ess}\cdot \sup(E(R, \alpha))$ for $\alpha > 0$ with  $R=\{\bd{X}_1, \cdots,\bd{X}_n\},$ where $\mbox{ess}\cdot \inf(Z)$ and  $\mbox{ess}\cdot \sup(Z)$ are the essential infinimum and the essential maximum of random variable $Z$, respectively.

For the  Tammes problem ($\alpha= -\infty$), the extremal energy  $\epsilon(R, -\infty)$ can be further studied through the random variable $\Theta_{\max}$.
Note that $\|\bd{x}_i-\bd{x}_j\|^2=2(1-\cos \theta_{ij}),$ where $\theta_{ij}$ is the angle between vectors $\overset{\longrightarrow}{\bd{O}\bd{x}_i}$ and $\overset{\longrightarrow}{\bd{O}\bd{x}_j}.$ Then
\beaa
\frac{1}{2 E(R, -\infty)^2}=\max_{\bd{x}_1, \cdots,\bd{x}_n \in \mathbb{S}^{p-1}}(1-\cos \theta_{ij})=1-\cos \tilde{\Theta}_{\max},
\eeaa
where $\tilde{\Theta}_{\max}=\max\{\theta_{ij};\, 1\leq i < j\leq n\}.$ Again, let $\bd{X}_1, \cdots, \bd{X}_n$ be i.i.d. random vectors with the uniform distribution on $\mathbb{S}^{p-1}$. Then, it is not difficult to see
\beaa
\frac{1}{2\epsilon(R, -\infty)^2}=\sup_{R}\frac{1}{2E(R, -\infty)^2}=\sup_{R}(1-\cos \tilde{\Theta}_{\max})=1-\cos \Delta
\eeaa
where $\Delta:=\mbox{ess}\cdot \sup(\Theta_{\max})$ is the essential upper bound of the random variable $\Theta_{\max}$ as defined in (\ref{fallway}).
Thus,
\bea\lbl{rhythm}
\epsilon(R, -\infty) = \frac{1}{\sqrt{2(1-\cos \Delta)}}.
\eea
The essential upper bound $\Delta$ of the random variable $\Theta_{\max}$ can be approximated by random sampling of $\Theta_{\max}$. So  the approach outlined above provides a direct way for using a stochastic method to study these deterministic problems and establishes connections between the random angles  and open problems mentioned above.
See, e.g.,  Katanforoush  and Shahshahani (2003) for further comments on randomization.   Recently, Armentano et al (2011) studied this problem by taking $\bd{x}_i$'s to be the roots of a special type of random polynomials.  Taking independent and uniform samples  $\bd{X}_1, \cdots, \bd{X}_n$ from the unit sphere $\mathbb{S}^{p-1}$ to get (\ref{rhythm}) is  simpler than using the roots of a random polynomials.

\section{Proofs}
\lbl{proof.sec}

\subsection{Technical Results}\lbl{suite'}

Recall that $\bd{X}_1, \bd{X}_2, \cdots$ are random points independently chosen with the uniform distribution on $\mathbb{S}^{p-1},$ the unit sphere in $\mathbb{R}^{p},$ and $\Theta_{ij}$  is the angle between $\overset{\longrightarrow}{\bd{O}\bd{X}_i}$ and $\overset{\longrightarrow}{\bd{O}\bd{X}_j}$ and $\rho_{ij}=\cos\Theta_{ij}$ for any $i\ne j.$ Of course, $\Theta_{ij}\in [0, \pi]$ for all $i \ne j.$ It is known that the distribution of $(\bd{X}_1, \bd{X}_2, \cdots)$ is the same as that of
\beaa
\Big(\frac{\bd{Y}_1}{\|\bd{Y}_1\|}, \frac{\bd{Y}_2}{\|\bd{Y}_2\|}, \cdots\Big)
\eeaa
where $\{\bd{Y}_1, \bd{Y}_2, \cdots\}$ are  independent $p$-dimensional random vectors with the normal distribution $N_p(\bd{0}, \bd{I}_p),$ that is, the normal distribution with mean vector $\bd{0}$ and the covariance matrix equal to  the $p\times p$ identity matrix $\bd{I}_p.$ Thus,
\beaa
\rho_{ij}=\cos\Theta_{ij}=\frac{\bd{Y}_i^T\bd{Y}_j}{\|\bd{Y}_i\|\cdot \|\bd{Y}_i\|}
\eeaa
for all $1\leq i<j \leq n.$ See, e.g., the Discussions in Section 5 from Cai and Jiang (2012) for further details. Of course, $\rho_{ii}=1$ and $|\rho_{ij}| \leq 1$ for all $i,j.$ Set
\bea\lbl{like}
M_n=\max_{1\leq i <  j \leq n}\rho_{ij}=\cos \Theta_{\min}.
\eea
\begin{lemma}\lbl{sentence}((22) in Lemma 4.2 from Cai and Jiang (2012))
Let $p\geq 2$. Then $\{\rho_{ij};\, 1\leq i < j\leq n\}$ are pairwise independent and identically distributed with  density function
\bea\lbl{walks}
g(\rho) =\frac{1}{\sqrt{\pi}}\frac{\Gamma(\frac{p}{2})}{\Gamma(\frac{p-1}{2})}\cdot(1-\rho^2)^{\frac{p-3}{2}},\ \ \ |\rho|<1.
\eea
\end{lemma}

Notice $y=\cos x$ is a strictly decresing function on $[0, \pi],$ hence $\Theta_{ij}=\cos^{-1}\rho_{ij}.$ A direct computation shows that Lemma \ref{sentence} is equivalent to the following lemma.

\begin{lemma}\lbl{recursion}
Let $p\geq 2$. Then,

(i) $\{\Theta_{ij};\, 1\leq i < j\leq n\}$ are pairwise independent and identically distributed with  density function
\bea\lbl{volume}
h(\theta) =\frac{1}{\sqrt{\pi}}\frac{\Gamma(\frac{p}{2})}{\Gamma(\frac{p-1}{2})}\cdot(\sin \theta)^{p-2},\ \ \ \theta \in [0, \pi].
\eea

(ii) If ``$\Theta_{ij}$" in (i) is replaced by ``$\pi -\Theta_{ij}$", the conclusion in (i) still holds.
\end{lemma}


Let $I$ be a finite set, and for each $\alpha\in I$, $X_{\alpha}$ be a Bernoulli random variable with $p_{\alpha}=P(X_{\alpha}=1)=1-P(X_{\alpha}=0)>0.$ Set $W=\sum_{\alpha\in I}X_{\alpha}$ and $\lambda=EW=\sum_{\alpha \in I}p_{\alpha}.$ For each $\alpha \in I,$ suppose we have chosen $B_{\alpha} \subset I$ with $\alpha \in B_{\alpha}.$ Define
\bea
b_1=\sum_{\alpha\in I}\sum_{\beta\in B_{\alpha}}p_{\alpha}p_{\beta}\ \ \mbox{and}\ \ b_2=\sum_{\alpha\in I}\sum_{\alpha \ne \beta\in B_{\alpha}}P(X_{\alpha}=1,\, X_{\beta}=1).
\eea

\begin{lemma}\lbl{speed}(Theorem 1 from Arratia et al. (1989)) For each $\alpha \in I,$ assume $X_{\alpha}$ is independent of $\{X_{\beta};\, \beta \in I-B_{\alpha}\}.$ Then $\big|P(X_{\alpha}=0\ \mbox{for all}\ \alpha\in I)-e^{-\lambda}\big|\leq b_1 + b_2.$
\end{lemma}

The following is essentially a special case of Lemma \ref{speed}.
\begin{lemma}\label{stein} Let $I$ be an index set and $\{B_{\alpha}, \alpha\in I\}$ be a set of subsets of $I,$ that is, $B_{\alpha}\subset I$ for each $\alpha \in I.$  Let also $\{\eta_{\alpha}, \alpha\in I\}$ be random variables. For a given $t\in \mathbb{R},$ set $\lambda=\sum_{\alpha\in I}P(\eta_{\alpha}>t).$ Then
\beaa
|P(\max_{\alpha \in I}\eta_{\alpha} \leq t)-e^{-\lambda}| \leq (1\wedge \lambda^{-1})(b_1+b_2+b_3)
\eeaa
where
\beaa
& & b_1=\sum_{\alpha \in I}\sum_{\beta \in B_{\alpha}}P(\eta_{\alpha} >t)P(\eta_{\beta} >t),\ \
 b_2=\sum_{\alpha \in I}\sum_{\alpha\ne \beta \in B_{\alpha}}P(\eta_{\alpha} >t, \eta_{\beta} >t),\\
 & & b_3=\sum_{\alpha \in I}E|P(\eta_{\alpha} >t|\sigma(\eta_{\beta}, \beta \notin B_{\alpha})) - P(\eta_{\alpha} >t)|,
\eeaa
and $\sigma(\eta_{\beta}, \beta \notin B_{\alpha})$ is the $\sigma$-algebra generated by $\{\eta_{\beta}, \beta \notin B_{\alpha}\}.$
In particular, if $\eta_{\alpha}$ is independent of $\{\eta_{\beta}, \beta \notin B_{\alpha}\}$ for each $\alpha,$ then $b_3=0.$
\end{lemma}

\begin{lemma}\lbl{newchildish} Let $p=p_n\geq 2$. Recall $M_n$ as in (\ref{like}). For $\{t_n\in [0,1];\, n\geq 2\}$, set
\beaa
h_n = \frac{n^2p^{1/2}}{\sqrt{2\pi}}
\int_{t_n}^1(1-x^2)^{\frac{p-3}{2}}\,dx.
\eeaa
If $\lim_{n\to\infty}p_n=\infty$ and $\lim_{n\to\infty}h_n=\lambda\in [0, \infty),$ then $\lim_{n\to\infty}P(M_n\leq t_n) = e^{-\lambda/2}.$
\end{lemma}
\textbf{Proof}. For brevity of notation, we sometimes write $t=t_n$ if there is no confusion. First, take $I=\{(i,j);\ 1\leq i< j \leq n\}.$ For $u =(i,j) \in I,$ set $B_{u}=\{(k,l)\in I;\ \mbox{one of}\ k\ \mbox{and}\ l =i\ \mbox{or}\ j,\ \mbox{but}\ (k,l)\ne u\},\ \eta_{u}=\rho_{ij}$ and $A_{u}=A_{ij}=\{\rho_{ij} > t\}.$  By the i.i.d. assumption on $\bd{X}_1, \cdots, \bd{X}_n$ and Lemma \ref{stein},
\bea\lbl{season}
|P(M_n\leq t) - e^{-\lambda_n}| \leq b_{1,n}+b_{2,n}
\eea
 where
\bea\lbl{sheep}
\lambda_n=\frac{n(n-1)}{2}P(A_{12})
\eea
 and
\beaa
b_{1,n}\leq 2n^3P(A_{12})^2\ \mbox{and}\ b_{2,n} \leq 2n^3P(A_{12}A_{13}).
\eeaa
By Lemma \ref{sentence}, $A_{12}$ and $A_{13}$ are independent events with the same probability. Thus, from (\ref{sheep}),
\bea\lbl{English}
b_{1,n} \vee b_{2,n} \leq 2n^3P(A_{12})^2 \leq \frac{8n\lambda_n^2}{(n-1)^2}\leq \frac{32\lambda_n^2}{n}
\eea
for all $n\geq 2.$ Now we compute $P(A_{12}).$ In fact, by Lemma \ref{sentence} again,
\beaa
P(A_{12})=\int_{t}^1g(x)\,dx &= & \frac{1}{\sqrt{\pi}}\frac{\Gamma(\frac{p}{2})}{\Gamma(\frac{p-1}{2})}
\int_{t}^1(1-x^2)^{\frac{p-3}{2}}\,dx.
\eeaa
Recalling the Stirling formula (see, e.g., p.368 from  Gamelin (2001) or (37) on
p.204 from Ahlfors (1979)):
\begin{eqnarray*}
\log\Gamma(z)=z\log z - z -\frac{1}{2}\log z+ \log \sqrt{2\pi}
 +O\left(\frac{1}{x}\right)
\end{eqnarray*}
as $x=\mbox{Re}\,(z)\to \infty,$ it is easy to verify that
\bea\lbl{white}
\frac{\Gamma(\frac{p}{2})}{\Gamma(\frac{p-1}{2})} \sim \sqrt{\frac{p}{2}}
\eea
as $p\to\infty.$ Thus,
\beaa
P(A_{12}) \sim \frac{p^{1/2}}{\sqrt{2\pi}}
\int_{t}^1(1-x^2)^{\frac{p-3}{2}}\,dx
\eeaa
as $n\to\infty.$ From (\ref{sheep}), we know
\beaa
\lambda_n \sim \frac{p^{1/2}n^2}{2\sqrt{2\pi}}
\int_{t}^1(1-x^2)^{\frac{p-3}{2}}\,dx=\frac{h_n}{2}
\eeaa
as $n\to\infty.$  Finally, by (\ref{season}) and (\ref{English}), we know
\beaa
\lim_{n\to\infty}P(M_n\leq t) = e^{-\lambda/2}\ \ \mbox{if}\ \ \lim_{n\to\infty}h_n=\lambda \in [0, \infty).\ \ \ \ \ \ \blacksquare
\eeaa

\subsection{Proofs of Main Results in Section \ref{fixedn}}\lbl{Proof.fixedn}

\begin{lemma}\lbl{kiss} Let $\bd{X}_1, \bd{X}_2,\cdots$ be independent random points with the uniform distribution on the unit sphere in $\mathbb{R}^{p}.$

(i) Let $p$ be fixed and $\mu$ be the probability measure with the density $h(\theta)$ as in (\ref{reliable}). Then, with probability one,
 $\mu_n$ in (\ref{empirical}) converges weakly to $\mu$ as $n\to\infty.$

(ii) Let $p=p_n$ and $\{\varphi_n(\theta);\, n\geq 1\}$ be a sequence of functions defined on $[0, \pi].$ If $\varphi_n(\Theta_{12})$ converges weakly to a probability measure $\nu$ as $n\to\infty,$ then, with probability one,
\bea\lbl{festival}
\nu_n:=\frac{1}{\binom{n}{2}}\sum_{1\leq i < j\leq n}\delta_{\varphi_n(\Theta_{ij})}
\eea
converges weakly to $\nu$ as $n\to\infty.$
\end{lemma}
\textbf{Proof}. First, we claim that, for any bounded and continuous function $u(x)$ defined on $\mathbb{R},$
\bea\lbl{TV}
\frac{1}{\binom{n}{2}}\sum_{1\leq i < j\leq n} \left[u(\varphi_n(\Theta_{ij})) -E u(\varphi_n(\Theta_{ij}))\right] \to 0\  \ a.s.
\eea
as $n\to\infty$ regardless $p$ is fixed  as in (i) or $p=p_n$ as in (ii) in the statement of the lemma. For convenience, write $u_n(\theta)=u(\varphi_n(\theta)).$ Then $u_n(\theta)$ is a bounded function with $M:=\sup_{\theta\in [0, \pi]}|u_n(\theta)|< \infty.$  By the Markov inequality
\begin{eqnarray*}
& & P\Big(\Big|\sum_{1\leq i < j\leq n} (u_n(\Theta_{ij})-Eu_n(\Theta_{ij}))
\Big| \geq \epsilon \binom{n}{2}\Big) \\
&\leq & \frac{1}{\binom{n}{2}^2\epsilon^2} E\Big|\sum_{1\leq i < j\leq n} (u_n(\Theta_{ij})-Eu_n(\Theta_{ij}))
\Big|^2
\end{eqnarray*}
for any $\epsilon>0.$ From (i) of  Lemma \ref{recursion},  $\{\Theta_{ij};\, 1\leq i < j\leq n\}$ are pairwise independent with the common distribution, the last expectation is therefore equal to $\binom{n}{2}\mbox{Var}(u_n(\Theta_{12})) \leq \binom{n}{2}M^2.$ This says that, for any $\epsilon>0,$
\[
P\Big(\Big|\sum_{1\leq i < j\leq n} (u_n(\Theta_{ij})-Eu_n(\Theta_{ij}))
\Big| \geq \epsilon \binom{n}{2}\Big) =O\Big(\frac{1}{n^2}\Big)
\]
as $n\to\infty.$ Note that the sum of the right hand side over all $n\geq 2$ is finite. By the Borel-Cantelli lemma, we conclude (\ref{TV}).

(i) Take $\varphi_n(\theta)=\theta$ for $\theta\in \mathbb{R}$ in (\ref{TV}) to get that
\begin{equation}\label{lions}
\frac{1}{\binom{n}{2}}\sum_{1\leq i < j\leq n} u(\Theta_{ij}) \to Eu(\Theta_{12})=\int_0^\pi u(\theta)h(\theta)\,d\theta\ \ a.s.
\end{equation}
as $n\to\infty,$ where $u(\theta)$ is any bounded continuous function on $[0, \pi]$ and $h(\theta)$ is as in (\ref{reliable}). This leads to that,  with probability one, $\mu_n$ in (\ref{empirical}) converges weakly to $\mu$ as $n\to\infty.$

(ii) Since  $\varphi_n(\Theta_{12})$ converges weakly to $\nu$ as $n\to\infty,$ we know that, for any bounded continuous function $u(x)$ defined on $\mathbb{R}$, $E u(\varphi_n(\Theta_{12})) \to \int_{-\infty}^{\infty}u(x)\,d\nu(x)$ as $n\to\infty.$ By (i) of Lemma \ref{recursion}, $E u(\varphi_n(\Theta_{ij}))=E u(\varphi_n(\Theta_{12}))$ for all $1\leq i < j\leq n.$ This and (\ref{TV}) yield
\begin{equation*}
\frac{1}{\binom{n}{2}}\sum_{1\leq i < j\leq n} u(\varphi_n(\Theta_{ij})) \to \int_{-\infty}^{\infty}u(x)\,d\nu(x)\ \ a.s.
\end{equation*}
as $n \to\infty$. Reviewing the definition of $\nu_n$ in (\ref{festival}), the above asserts that, with probability one, $\nu_n$ converges weakly to $\nu$ as $n \to\infty.$\ \ \ \ \ \ \  $\blacksquare$\\

\noindent\textbf{Proof of Theorem \ref{quarter}}. This is a direct consequence of (i) of Lemma \ref{kiss}.   \ \ \ \ \ \ \ \ $\blacksquare$\\

Recall $\bd{X}_1, \cdots, \bd{X}_n$ are random points independently chosen with the uniform distribution on $\mathbb{S}^{p-1},$ the unit sphere in $\mathbb{R}^{p},$ and $\Theta_{ij}$  is the angle between $\overset{\longrightarrow}{\bd{O}\bd{X}_i}$ and $\overset{\longrightarrow}{\bd{O}\bd{X}_j}$ and $\rho_{ij}=\cos\Theta_{ij}$ for all $1\leq i, j\leq n.$ Of course, $\rho_{ii}=1$ and $|\rho_{ij}| \leq 1$ for all $1\leq i \ne j \leq n.$ Review (\ref{like}) to have
\beaa
M_n=\max_{1\leq i <  j \leq n}\rho_{ij}=\cos \Theta_{\min}.
\eeaa
To prove Theorem \ref{tube}, we need the following result.
\begin{prop}\lbl{probability} Fix $p\geq 2.$ Then $n^{4/(p-1)}(1-M_n)$ converges  to the distribution function
\[\lbl{chick2}
F_1(x)=1-\exp\{-K_1x^{(p-1)/2}\},\ \ x\geq 0,
\]
in distribution as $n\to \infty,$ where
\bea\lbl{luck1}
K_1=\frac{2^{(p-5)/2}}{\sqrt{\pi}}\frac{\Gamma(\frac{p}{2})}{\Gamma(\frac{p+1}{2})}.
\eea
\end{prop}

\noindent\textbf{Proof}. Set $t=t_n=1-xn^{-4/(p-1)}$ for $x\geq 0.$ Then
\bea\lbl{fur}
 t\to 1\ \ \mbox{and}\ \
 t^2 =1 - \frac{2x}{n^{4/(p-1)}} + O\Big(\frac{1}{n^{8/(p-1)}}\Big)
\eea
as $n\to\infty.$ Notice
\beaa
P(n^{4/(p-1)}(1-M_n) < x)=P(M_n> t)=1-P(M_n \leq t).
\eeaa
Thus, to prove the theorem, since $F_1(x)$ is continuous, it is enough to show that
\bea\lbl{hate}
P(M_n \leq t) \to e^{-K_1 x^{(p-1)/2}}
\eea
as $n\to\infty,$ where $K_1$ is as in (\ref{luck1}).

Now, take $I=\{(i,j);\ 1\leq i< j \leq n\}.$ For $u =(i,j) \in I,$ set $B_{u}=\{(k,l)\in I;\ \mbox{one of}\ k\ \mbox{and}\ l =i\ \mbox{or}\ j,\ \mbox{but}\ (k,l)\ne u\},\ \eta_{u}=\rho_{ij}$ and $A_{u}=A_{ij}=\{\rho_{ij} > t\}.$  By the i.i.d. assumption on $\bd{X}_1, \cdots, \bd{X}_n$ and Lemma \ref{stein},
\bea\lbl{season1}
|P(M_n\leq t) - e^{-\lambda_n}| \leq b_{1,n}+b_{2,n}
\eea
 where
\bea\lbl{sheep2}
\lambda_n=\frac{n(n-1)}{2}P(A_{12})
\eea
 and
\beaa
b_{1,n}\leq 2n^3P(A_{12})^2\ \mbox{and}\ b_{2,n} \leq 2n^3P(A_{12}A_{13}).
\eeaa
By Lemma \ref{sentence}, $A_{12}$ and $A_{13}$ are independent events with the same probability. Thus, from (\ref{sheep2}),
\bea\lbl{English1}
b_{1,n} \vee b_{2,n} \leq 2n^3P(A_{12})^2 \leq \frac{8n\lambda_n^2}{(n-1)^2}\leq \frac{32\lambda_n^2}{n}
\eea
for all $n\geq 2.$ Now we evaluate $P(A_{12}).$ In fact, by Lemma \ref{sentence} again,
\beaa
P(A_{12})=\int_{t}^{1}g(x)\,dx &= & \frac{1}{\sqrt{\pi}}\frac{\Gamma(\frac{p}{2})}{\Gamma(\frac{p-1}{2})}
\int_{t}^1(1-x^2)^{\frac{p-3}{2}}\,dx.
\eeaa
Set $m=\frac{p-3}{2}\geq -\frac{1}{2}.$ We claim
\bea\lbl{way1}
\int_{t}^1(1-x^2)^{m}\,dx \sim \frac{1}{2m+2}(1-t^2)^{m+1}
\eea
as $n\to\infty.$ In fact,  set $s=x^2.$ Then $x=\sqrt{s}$ and $dx=\frac{1}{2\sqrt{s}} ds.$ It follows that
\beaa
\int_{t}^1(1-x^2)^{m}\,dx & = & \int_{t^2}^1 \frac{1}{2\sqrt{s}}(1-s)^{m}\,ds \\
 & \sim & \frac{1}{2}\int_{t^2}^1 (1-s)^{m}\,ds=\frac{1}{2m+2}(1-t^2)^{m+1}
\eeaa
as $n\to\infty,$ where the fact $\lim_{n\to\infty} t=\lim_{n\to\infty} t_n=1$ stated in (\ref{fur}) is used in the second step to replace $\frac{1}{2\sqrt{s}}$ by $\frac{1}{2}.$ So the claim (\ref{way1}) follows.

Now, we know from (\ref{sheep2}) that
\beaa
\lambda_n \sim \frac{n^2}{2\sqrt{\pi}}\frac{\Gamma(\frac{p}{2})}{\Gamma(\frac{p-1}{2})}
\int_{t}^1(1-x^2)^{\frac{p-3}{2}}\,dx
& \sim & \frac{n^2}{2\sqrt{\pi}}\frac{\Gamma(\frac{p}{2})}{(p-1)\Gamma(\frac{p-1}{2})}(1-t^2)^{(p-1)/2}\\
& = & \frac{1}{4\sqrt{\pi}}\frac{\Gamma(\frac{p}{2})}{\Gamma(\frac{p+1}{2})}\left(n^{4/(p-1)}(1-t^2)\right)^{(p-1)/2}
\eeaa
as $n\to\infty,$  where (\ref{way1}) is used in the second step and the fact $\Gamma(x+1)=x\Gamma(x)$ is used in the last step. By (\ref{fur}),
\[
n^{4/(p-1)}(1-t^2)=2x + O\Big(\frac{1}{n^{4/(p-1)}}\Big)
\]
as $n\to\infty.$ Therefore,
\[
\lambda_n \to \frac{2^{(p-5)/2}}{\sqrt{\pi}}\frac{\Gamma(\frac{p}{2})}{\Gamma(\frac{p+1}{2})} x^{(p-1)/2}=K_1 x^{(p-1)/2}
\]
as $n\to\infty.$ Finally, by (\ref{season1}) and (\ref{English1}), we know
\beaa
\lim_{n\to\infty}P(M_n\leq t) = e^{-K_1 x^{(p-1)/2}}.
\eeaa
 This concludes (\ref{hate}).\ \ \ \ \ \ \ \ $\blacksquare$\\

\noindent\textbf{Proof of Theorem \ref{tube}}. First, since $M_n=\cos \Theta_{\min}$ by (\ref{geese}), then use the identity $1-\cos h =2\sin^2\frac{h}{2}$ for all $h\in \mathbb{R}$ to have
\bea\lbl{kick}
n^{4/(p-1)}(1-M_n)=2n^{4/(p-1)}\sin^2\frac{\Theta_{\min}}{2}.
\eea
By Proposition \ref{probability} and the Slusky lemma, $\sin\frac{\Theta_{\min}}{2}\to 0$ in probability as $n\to\infty.$  Noticing $0\leq \Theta_{\min}\leq \pi$, we then have $\Theta_{\min} \to 0$ in probability as $n\to\infty.$ From (\ref{kick}) and the fact that $\lim_{x\to 0}\frac{\sin x}{x}=1$ we obtain
\beaa
\frac{n^{4/(p-1)}(1-M_n)}{\frac{1}{2} n^{4/(p-1)}\Theta_{\min}^2} \to 1
\eeaa
in probability as $n\to\infty.$ By Proposition \ref{probability} and the Slusky lemma again, $\frac{1}{2}n^{4/(p-1)}\Theta_{\min}^2$ converges in distribution to $F_1(x)$ as in Proposition \ref{probability}. Second, for any $x>0,$
\bea\lbl{return}
P(n^{2/(p-1)}\Theta_{\min} \leq x)
& = & P\Big(\frac{1}{2}n^{4/(p-1)}\Theta_{\min}^2 \leq \frac{x^2}{2}\Big)\nonumber\\
& \to & 1-\exp\{-K_1(x^2/2)^{(p-1)/2}\}=1-\exp\{-Kx^{p-1}\}
\eea
as $n\to\infty,$ where
\bea\lbl{road}
K=2^{(1-p)/2}K_1=\frac{1}{4\sqrt{\pi}}\frac{\Gamma(\frac{p}{2})}{\Gamma(\frac{p+1}{2})}.
\eea
Now we prove
\bea\lbl{scholar}
n^{2/(p-1)}(\pi-\Theta_{\max})\ \mbox{converges weakly to}\ F(x)\ \mbox{as}\ n\to\infty.
\eea
 In fact, recalling the proof of the above and that of Proposition \ref{probability}, we only use the following  properties about $\rho_{ij}:$

(i) $\{\rho_{ij};\, 1\leq i < j\leq n\}$ are pairwise independent.

(ii) $\rho_{ij}$ has density function $g(\rho)$ given in (\ref{walks}) for all $1\leq i< j\leq n.$

(iii) For each $1\leq i< j\leq n$,  $\rho_{ij}$ is independent of $\{\rho_{kl};\, 1\leq k<l\leq n;\, \{k, l\} \cap \{i,j\}=\emptyset\}.$

By using Lemmas \ref{sentence} and \ref{recursion} and the remark between them, we see that the above three properties are equivalent to

(a) $\{\Theta_{ij};\, 1\leq i < j\leq n\}$ are pairwise independent.

(b) $\Theta_{ij}$ has density function $h(\theta)$ given in (\ref{volume}) for all $1\leq i< j\leq n.$

(c) For each $1\leq i< j\leq n$,  $\Theta_{ij}$ is independent of $\{\Theta_{kl};\, 1\leq k<l\leq n;\, \{k, l\} \cap \{i,j\}=\emptyset\}.$

It is easy to see from (ii) of lemma \ref{recursion} that the above three properties are equivalent to the corresponding (a) , (b) and (c) when ``$\Theta_{ij}$" is replaced by ``$\pi-\Theta_{ij}$" and ``$\Theta_{kl}$" is replaced by ``$\pi-\Theta_{kl}.$" Also, it is key to observe that $\min\{\pi-\Theta_{ij};\, 1\leq i < j\leq n\}=\pi - \Theta_{\max}.$ We then deduce from (\ref{return}) that
\bea\lbl{return1}
P(n^{2/(p-1)}(\pi - \Theta_{\max}) \leq x)\to 1-\exp\{-Kx^{p-1}\}
\eea
as $n\to\infty,$ where $K$ is as in (\ref{road}).\ \ \ \ \ \ \ $\blacksquare$\\

\noindent\textbf{Proof of Theorem \ref{goose}}. We will prove the following:
\bea\lbl{Nash}
\lim_{n\to\infty}P\big(n^{2/(p-1)}\Theta_{\min}\geq x,\, n^{2/(p-1)}(\pi-\Theta_{\max}) \geq y\big)= e^{-K(x^{p-1}+y^{p-1})}
\eea
for any $x\geq 0$ and $y\geq 0,$ where $K$ is as in (\ref{luck}). Note that the right hand side in (\ref{Nash}) is identical to $P(X\geq x,\, Y \geq y),$ where $X$ and $Y$ are as in the statement of Theorem \ref{goose}. If (\ref{Nash}) holds, by the fact that $\Theta_{\min}, \Theta_{\max}, X, Y$ are continuous random variables and by Theorem \ref{tube} we know that   $Q_n:=\big((n^{2/(p-1)}\Theta_{\min}, n^{2/(p-1)}(\pi-\Theta_{\max})\big)\in \mathbb{R}^2$ for $n\geq 2$  is a tight sequence. By the standard subsequence argument, we obtain that $Q_n$ converges weakly to the distribution of $(X, Y)$ as $n\to\infty.$ Applying the map $h(x,y)=x-y$ with $x, y \in \mathbb{R}$ to the sequence $\{Q_n;\, n\geq 2\}$ and its limit, the desired conclusion then follows from the continuous mapping theorem on the weak convergence of probability measures.

We now prove (\ref{Nash}). Set $t_x=n^{-2/(p-1)}x\ \ \mbox{and}\ \ t_y=\pi-n^{-2/(p-1)}y.$  Without loss of generality, we assume $0\leq t_x<t_y<\infty$ for all $n\geq 2.$ Then
\bea\lbl{doubt}
& & P\big(n^{2/(p-1)}\Theta_{\min}\geq x,\, n^{2/(p-1)}(\pi-\Theta_{\max}) \geq y\big)\nonumber\\
&= & P(t_x\leq \Theta_{ij}\leq t_y\ \mbox{for all}\, 1\leq i <j \leq n)\nonumber\\
&= & P\big(X_u=0\ \mbox{for all}\ u\in I\big)
\eea
where $I:=\{(i,j);\ 1\leq i< j \leq n\}$ and
\beaa
X_{u}:=\begin{cases}
1, & \text{if $\Theta_u\notin [t_x, t_y]$};\\
0, & \text{if $\Theta_u\in [t_x, t_y].$}
\end{cases}
\eeaa
For $u =(i,j) \in I,$ set $B_{u}=\{(k,l)\in I;\ \mbox{one of}\ k\ \mbox{and}\ l =i\ \mbox{or}\ j,\ \mbox{but}\ (k,l)\ne u\}.$   By the i.i.d. assumption on $\bd{X}_1, \cdots, \bd{X}_n$ and Lemma \ref{speed}
\bea\lbl{far}
|P\big(X_u=0\ \mbox{for all}\ u\in I\big) - e^{-\lambda_n}| \leq b_{1,n}+b_{2,n}
\eea
 where
\bea\lbl{sheep5}
\lambda_n=\frac{n(n-1)}{2}P(A_{12})\ \ \mbox{and}\ \ A_{12}=\big\{\Theta_{12}\notin [t_x, t_y] \big\}
\eea
 and
\bea\lbl{california}
b_{1,n}\leq 2n^3P(A_{12})^2\ \mbox{and}\ b_{2,n} \leq 2n^3P(A_{12}A_{13})=2n^3P(A_{12})^2
\eea
by Lemma \ref{recursion}. Now
\bea\lbl{addition}
P(A_{12})=P(\Theta_{12}< t_x) + P(\Theta_{12}> t_y).
\eea
By Lemma \ref{recursion} again,
\bea
P(\Theta_{12}> t_y) &= & \frac{1}{\sqrt{\pi}}\frac{\Gamma(\frac{p}{2})}{\Gamma(\frac{p-1}{2})}
\int_{t_y}^{\pi}(\sin \theta)^{p-2}\,d\theta\nonumber\\
& = & \frac{1}{\sqrt{\pi}}\frac{\Gamma(\frac{p}{2})}{\Gamma(\frac{p-1}{2})}
\int^{n^{-2/(p-1)}y}_{0}(\sin \eta)^{p-2}\,d\eta\lbl{realize}
\eea
by setting $\eta=\pi -\theta.$ Now, set $v=\cos \eta$ for $\eta\in [0, \pi].$ Write $(\sin \eta)^{p-2}=-(\sin \eta)^{p-3}(\cos \eta)'.$ Then  the  integral in (\ref{realize}) is equal to
\beaa
\int_{v_y}^{1}(1-v^2)^{(p-3)/2}\,dv
\eeaa
where
\beaa
v_y:=\cos (n^{-2/(p-1)}y)=1-\frac{y^2}{2n^{4/(p-1)}} + O\Big(\frac{1}{n^{8/(p-1)}}\Big)
\eeaa
as $n\to\infty$ by the Taylor expansion. Trivially,
\beaa\lbl{war}
v_y^2=1-\frac{y^2}{n^{4/(p-1)}} + O\Big(\frac{1}{n^{8/(p-1)}}\Big)
\eeaa
as $n\to\infty.$ Thus, by (\ref{way1}),
\beaa
\int_{v_y}^{1}(1-v^2)^{(p-3)/2}\,dv \sim \frac{1}{p-1}(1-v_y^2)^{(p-1)/2}=\frac{y^{p-1}}{(p-1)n^{2}}\Big(1+ O\Big(\frac{1}{n^{4/(p-1)}}\Big)\Big)
\eeaa
as $n\to\infty.$ Combining all the above we conclude that
\bea\lbl{drop}
P(\Theta_{12}> t_y) & = & \frac{\Gamma(\frac{p}{2})}{\sqrt{\pi}(p-1)\Gamma(\frac{p-1}{2})}\frac{y^{p-1}}{n^{2}}(1+o(1))\nonumber\\
& = & \frac{\Gamma(\frac{p}{2})}{2\sqrt{\pi}\,\Gamma(\frac{p+1}{2})}\frac{y^{p-1}}{n^{2}}(1+o(1))
\eea
as $n\to\infty.$ Similar to the part between (\ref{realize}) and (\ref{drop}), we have
\beaa
P(\Theta_{12}< t_x) &= & \frac{1}{\sqrt{\pi}}\frac{\Gamma(\frac{p}{2})}{\Gamma(\frac{p-1}{2})}
\int^{n^{-2/(p-1)}x}_{0}(\sin \theta)^{p-2}\,d\theta\\
& = & \frac{\Gamma(\frac{p}{2})}{2\sqrt{\pi}\,\Gamma(\frac{p+1}{2})}\frac{x^{p-1}}{n^{2}}(1+o(1))
\eeaa
as $n\to\infty.$ This joint with (\ref{drop}) and (\ref{addition}) implies that
\beaa
P(A_{12})=\frac{\Gamma(\frac{p}{2})}{2\sqrt{\pi}\,\Gamma(\frac{p+1}{2})}\frac{x^{p-1}+y^{p-1}}{n^{2}}(1+o(1))
\eeaa
as $n\to\infty.$ Recalling (\ref{sheep5}) and (\ref{california}), we obtain
\beaa
\lim_{n\to\infty}\lambda_n=K(x^{p-1}+y^{p-1})
\eeaa
and $b_{1,n} \vee b_{2,n} =O\Big(\frac{1}{n}\Big)$ as $n\to\infty,$ where $K$ is as in (\ref{luck}). These two assertions and (\ref{far}) yield
\beaa
\lim_{n\to\infty}P\big(X_u=0\ \mbox{for all}\ u\in I\big)=e^{-K(x^{p-1} + y^{p-1})}.
\eeaa
Finally, this together with (\ref{doubt}) implies (\ref{Nash}).\ \ \ \ \ \ $\blacksquare$

\subsection{Proofs of Main Results in Section \ref{largen}}\lbl{Proof.largen}

%
%
%

\noindent\textbf{Proof of Theorem \ref{nickle}}. Notice  $(p-2)/p\to 1$ as $p\to\infty,$ to prove the theorem, it is enough to show that the theorem holds if ``$\mu_{n,p}$" is replaced by ``$\frac{1}{\binom{n}{2}}\sum_{1\leq i < j\leq n}\delta_{\sqrt{p}({\frac{\pi}{2}-\Theta_{ij}})}.$" Thus,   without loss of generality, we assume (with a bit of abuse of notation) that
\bea\lbl{empirical5}
 \mu_{n,p}=\frac{1}{\binom{n}{2}}\sum_{1\leq i < j\leq n}\delta_{\sqrt{p}({\frac{\pi}{2}-\Theta_{ij}})},\ n \geq 2,\ p\geq 2.
 \eea
Recall $p=p_n.$ Set $Y_n:=\sqrt{p}(\frac{\pi}{2}-\Theta_{12})$ for $p\geq 2.$ We claim that
\bea\lbl{domain}
Y_n\ \mbox{converges weakly to}\ N(0,1)
\eea
as $n\to\infty.$ Assuming this is true, taking $\varphi_n(\theta)=\sqrt{p}(\frac{\pi}{2}-\theta)$ for $\theta \in [0, \pi]$ and $\nu=N(0,1)$ in (ii) of Lemma \ref{kiss}, then, with probability one, $\mu_{n,p}$  converges weakly to $N(0, 1)$ as $n\to \infty.$

Now we prove the claim. In fact, noticing $\Theta_{12}$ has density
$h(\theta)$ in (\ref{volume}), it is easy to see that $Y_n$ has density function
\bea
h_n(y): &= & \frac{1}{\sqrt{\pi}}\frac{\Gamma(\frac{p}{2})}{\Gamma(\frac{p-1}{2})}\cdot\Big[\sin \Big(\frac{\pi}{2}-\frac{y}{\sqrt{p}}\Big)\Big]^{p-2}\cdot \Big|-\frac{1}{\sqrt{p}}\Big| \nonumber\\
 & = & \frac{1}{\sqrt{p\pi}}\frac{\Gamma(\frac{p}{2})}{\Gamma(\frac{p-1}{2})}\cdot\Big(\cos \frac{y}{\sqrt{p}}\Big)^{p-2}\lbl{manager}
\eea
for any $y\in \mathbb{R}$ as $n$ is sufficiently large since $\lim_{n\to\infty}p_n=\infty.$ By (\ref{white}),
\bea\lbl{loom}
\frac{1}{\sqrt{p\pi}}\frac{\Gamma(\frac{p}{2})}{\Gamma(\frac{p-1}{2})} \to \frac{1}{\sqrt{2\pi}}
\eea
as $n\to\infty.$ On the other hand, by the Taylor expansion,
\beaa
\Big(\cos \frac{y}{\sqrt{p}}\Big)^{p-2}=\Big(1-\frac{y^2}{2p} + O\Big(\frac{1}{p^2}\Big)\Big)^{p-2} \to e^{-y^2/2}
\eeaa
as $n\to\infty.$ The above together with (\ref{manager}) and (\ref{loom}) yields that
\bea\lbl{tipping}
\lim_{n\to\infty} h_n(y) \to \frac{1}{\sqrt{2\pi}}e^{-y^2/2}
\eea
 for any $y\in \mathbb{R}.$ The assertions in  (\ref{manager}) and (\ref{loom}) also imply that $\sup_{y\in \mathbb{R}}|h_n(y)|\leq C$ for $n$ sufficiently large, where $C$ is a constant not depending on $n.$ This and (\ref{tipping}) conclude (\ref{domain}). \ \ \ \ \ \  $\blacksquare$\\

\noindent\textbf{Proof of Proposition \ref{rabbit}}. By (i) of Lemma \ref{recursion},
\beaa
P(|\Theta -\frac{\pi}{2}|\geq \epsilon)=C_p\int_{|\theta -\frac{\pi}{2}|\geq \epsilon}(\sin \theta)^{p-2}\,d\theta=C_p\int_{\epsilon\leq |t| \leq \pi/2} (\cos t)^{p-2}\,dt
\eeaa
by making transform $t=\theta-\frac{\pi}{2},$ where $C_p:=\frac{1}{\sqrt{\pi}}\Gamma(\frac{p}{2})/\Gamma(\frac{p-1}{2}).$ The last term above is identical to
\beaa
2C_p\int_{\epsilon}^{\pi/2}(\cos t)^{p-2}\,dt \leq \pi C_p(\cos \epsilon)^{p-2}.
\eeaa
It is known that $\lim_{x\to +\infty}\Gamma(x+a)/(x^a\Gamma(x))= 1,$ see, e.g., Dong, Jiang and Li (2012). Then $\pi C_p\leq K\sqrt{p}$ for all $p\geq 2,$ where $K$ is a universal constant. The desired conclusion then follows. \ \ \ \ \ \ \ $\blacksquare$\\

\noindent\textbf{Proof of Theorem \ref{newkind}}. Review the proof of Theorem 1 in Cai and Jiang (2012). Replacing $|\rho_{ij}|$, $L_n$ in (2) and Lemma 6.4  from Cai and Jiang (2012) with $\rho_{ij}$, $M_n$ in (\ref{like}) and Lemma \ref{newchildish} here, respectively. In the places where ``$n-2$" or ``$n-4$" appear in the proof, change them to ``$p-1$" or ``$p-3$" accordingly. Keeping the same argument in the proof, we then obtain the following.

 (a) $M_n\to 0$ in probability as $n\to\infty.$

 (b) Let $T_n=\log (1-M_n^2).$ Then, as $n\to\infty,$
\beaa
pT_n + 4\log n -\log \log n
\eeaa
 converges weakly to an extreme value distribution with the distribution function $F(y)=1- e^{-Ke^{y/2}},\  y\in\mathbb{R}$ and $K=1/(2\sqrt{8\pi})=1/(4\sqrt{2\pi}).$ From (\ref{like}) we know
\bea
& & M_n=\max_{1\leq i <  j \leq n}\rho_{ij}=\cos \Theta_{\min}\ \ \mbox{and}\ \ \Theta_{\min}\in [0, \pi];\lbl{chewing}\\
& & T_n=\log (1-M_n^2)=2\log \sin \Theta_{\min}.\lbl{lyrics}
\eea
Then (a) above implies that $\Theta_{\min} \to \pi/2$ in probability as $n\to\infty,$ and (b) implies (ii) for $ \Theta_{\min}$ in the statement of Theorem \ref{newkind}. Now, observe that
\bea\lbl{remote}
\min_{1\leq i< j \leq n}\{\pi-\Theta_{ij}\}=\pi- \Theta_{\max}\ \ \mbox{and}\ \ \sin(\pi-\Theta_{\max})=\sin \Theta_{\max}.
\eea
By the same argument between (\ref{scholar}) and (\ref{return1}), we get $\pi-\Theta_{\max} \to \pi/2$ in probability as $n\to\infty,$  that is, $\Theta_{\max} \to \pi/2$ in probability as $n\to\infty.$ Notice
\beaa
& & \max_{1\leq i < j\leq p}\Big|\Theta_{ij} - \frac{\pi}{2}\Big|\\
& \leq & \Big|\Theta_{\max} -\frac{\pi}{2}\Big|+  \Big|\Theta_{\min} -\frac{\pi}{2}\Big|\to 0
\eeaa
in probability as $n\to\infty.$ We get (i).

Finally, by the same argument between (\ref{scholar}) and (\ref{return1}) again, and by (\ref{remote}) we obtain
\beaa
2p\log \sin \Theta_{\max} + 4\log n -\log \log n
\eeaa
converges weakly to  $F(y)=1- e^{-Ke^{y/2}},\  y\in\mathbb{R}$ and $K=1/(4\sqrt{2\pi}\,).$ Thus, (ii) also holds for $ \Theta_{\max}.$\ \ \ \ \ \ $\blacksquare$\\

\noindent\textbf{Proof of Corollary \ref{newbrother}}. Review the proof of Corollary 2.2 from Cai and Jiang (2012). Replacing $L_n$ and Theorem 1 there by $M_n$ and Theorem \ref{newkind}, we get that
\beaa
pM_n^2 - 4\log n +\log \log n
\eeaa
converges weakly to the distribution function $\exp\{-{1\over 4\sqrt{2\pi}} e^{-(y+8\alpha^2)/2}\}$,  $y\in \mathbb{R}$. The desired conclusion follows since $M_n= \cos \Theta_{\min}.$\ \ \ \ \ \ \ $\blacksquare$\\

\noindent\textbf{Proof of Theorem \ref{newjean}}. Review the proof of Theorem 2 in Cai and Jiang (2012). Replacing $|\rho_{ij}|$, $L_n$ in (2) and Lemma 6.4  from Cai and Jiang (2012) with $\rho_{ij}$, $M_n$ in (\ref{like}) and Lemma \ref{newchildish}, respectively. In the places where ``$n-2$" and ``$n-4$" appear in the proof, change them to ``$p-1$" and ``$p-3$" accordingly. Keeping the same argument in the proof, we then have the following conclusions.

(i) $M_n\to \sqrt{1-e^{-4\beta}}$ in probability as $n\to\infty.$

(ii) Let $T_n=\log (1-M_n^2).$ Then, as $n\to\infty,$
\beaa
pT_n + 4\log n -\log \log n
\eeaa
 converges weakly to the distribution function
\bea
F(y)=1- \exp\left\{-K(\beta) e^{(y+8\beta)/2}\right\},\  y\in\mathbb{R},
\eea
where
\beaa
K(\beta)=\frac{1}{2}\Big(\frac{\beta}{2\pi(1- e^{-4\beta})}\Big)^{1/2}=\Big(\frac{\beta}{8\pi(1- e^{-4\beta})}\Big)^{1/2}.
\eeaa
From (\ref{chewing}) and (\ref{lyrics}) we obtain
\bea
& & \Theta_{\min}\to \cos^{-1}\sqrt{1-e^{-4\beta}}\ \mbox{in probability and}\lbl{payment}\\
& & 2p\log \sin \Theta_{\min} + 4\log n -\log \log n\lbl{graduation}
\eea
converges weakly to the distribution function
\bea\lbl{globe}
F(y)=1- \exp\left\{-K(\beta) e^{(y+8\beta)/2}\right\},\  y\in\mathbb{R},\ \mbox{where}\    K(\beta)=\Big(\frac{\beta}{8\pi(1- e^{-4\beta})}\Big)^{1/2}
\eea
as $n\to\infty$. Now, reviewing (\ref{remote}) and  the argument between (\ref{scholar}) and (\ref{return1}), by (\ref{payment}) and (\ref{graduation}), we conclude that
$\Theta_{\max}\to \pi-\cos^{-1}\sqrt{1-e^{-4\beta}}$ in probability  and $2p\log \sin \Theta_{\max} + 4\log n -\log \log n$ converges weakly to the distribution function $F(y)$ as in (\ref{globe}). The proof is completed. \ \ \ \ \ \ $\blacksquare$\\

\noindent\textbf{Proof of Theorem \ref{newlearn}}. Review the proof of Theorem 3 in Cai and Jiang (2012). Replacing $|\rho_{ij}|$, $L_n$ in (2) and Lemma 6.4  from Cai and Jiang (2012) with $\rho_{ij}$, $M_n$ in (\ref{like}) and Lemma \ref{newchildish}, respectively. In the places where ``$n-2$" or ``$n-4$" appear in the proof, change them to ``$p-1$" or ``$p-3$" accordingly. Keeping the same argument in the proof, we get the following results.

i) $M_n \to 1$ in probability as $n\to\infty.$

ii) As $n\to\infty,$
\beaa
pM_n + \frac{4p}{p-1}\log n-\log p
\eeaa
 converges weakly to the distribution function $F(y)=1- e^{-Ke^{y/2}},\  y\in\mathbb{R}$ with $K=1/(2\sqrt{2\pi}).$ Combining i), ii), (\ref{chewing}) and (\ref{lyrics}), we see that, as $n\to\infty$,
\beaa
& & \Theta_{\min}\to 0\  \mbox{in probability};\\
& & 2p\log \sin \Theta_{\min} + \frac{4p}{p-1}\log n-\log p\ \mbox{converges weakly to}\ \ \ \ \ \ \ \ \ \ \ \ \ \ \ \ \
\eeaa
$F(y)=1- e^{-Ke^{y/2}},\  y\in\mathbb{R}$ with $K=1/(2\sqrt{2\pi}).$ Finally, combining the above two convergence results with (\ref{remote}) and the argument between (\ref{scholar}) and (\ref{return1}), we have
\beaa
& & \Theta_{\max}\to \pi\ \mbox{in probability};\\
& & 2p\log \sin \Theta_{\max} + \frac{4p}{p-1}\log n-\log p\ \mbox{converges weakly to}\ \ \ \ \ \ \ \ \ \ \ \ \ \ \ \ \
\eeaa
$F(y)=1- e^{-Ke^{y/2}},\  y\in\mathbb{R}$ with $K=1/(2\sqrt{2\pi}).$\ \ \ \ \ \ $\blacksquare$\\


\end{document}